\documentclass[onecolumn
aps,
prf,
showkeys
]{revtex4-2}

\usepackage{amsmath,amsfonts}
\usepackage{tikz}
\usetikzlibrary{matrix}
\usepackage{import}
\usepackage{multirow}
\usepackage{float}

\usepackage{graphicx}
\usepackage{dcolumn}
\usepackage{bm}
\usepackage{hyperref}

\usepackage[noabbrev,nameinlink]{cleveref}

\usepackage[dvipsnames]{xcolor}
\definecolor{mycolor1}{RGB}{24, 92, 68}

\hypersetup{
    colorlinks = true,
    linkcolor = {mycolor1},
    citecolor ={mycolor1}
}

\usepackage{orcidlink} 

\newcolumntype{C}[1]{>{\centering\arraybackslash}m{#1}} 
\DeclareMathOperator*{\argmin}{argmin}

\begin{document}

\preprint{APS/123-QED}

\title{Parametric shape optimization of flagellated microswimmers using Bayesian techniques}

\author{Lucas Palazzolo \orcidlink{0009-0002-4188-014X}}
 \email{lucas.palazzolo@inria.fr}
\affiliation{Université Côte d'Azur, Inria, Calisto team, Sophia Antipolis, France}

\author{Mickaël Binois\orcidlink{0000-0002-7225-1680}}
\email{mickael.binois@inria.fr}
\affiliation{Université Côte d'Azur, Inria, Acumes team, CNRS, LJAD, Sophia Antipolis, France}

\author{Luca Berti\orcidlink{0000-0001-9149-2271}}
\affiliation{Cemosis, IRMA UMR 7501, CNRS, Université de Strasbourg, France}

\author{Laetitia Giraldi\orcidlink{0000-0003-2684-0203}}%
\email{laetitia.giraldi@inria.fr}
\affiliation{Université Côte d'Azur, Inria, Calisto team, Sophia Antipolis, France~}

\date{\today}

\begin{abstract}
Understanding and optimizing the design of helical microswimmers is crucial for advancing their application in various fields. This study presents an innovative approach combining Free-Form Deformation with Bayesian Optimization to enhance the shape of these swimmers. Our method facilitates the computation of generic swimmer shapes that achieve optimal average speed and efficiency. Applied to both monoflagellated and biflagellated swimmers, our optimization framework has led to the identification of new optimal shapes. These shapes are compared with biological counterparts, highlighting a diverse range of swimmers, including both \textit{pushers} and \textit{pullers}.
\end{abstract}


\maketitle


\section{Introduction}

Microswimmers are microscopic organisms with the ability to move in a fluid environment, as seen in sperm cells, bacteria, and other similar entities. Understanding and improving their performance, including their swimming mechanisms and shapes, is a current research focus. This knowledge can be applied to the design of microrobots with a wide range of applications, such as delivering drugs to specific target locations in medicine \cite{review_robots_medical_2019, medical_robots_cancer_2020}.\\

\noindent Because of their small size, microswimmers operate at low Reynolds numbers, within the Stokes regime. As explained by Purcell in the \textit{scallop theorem}, swimmers performing reciprocal motion do not produce any net displacement at low Reynolds numbers, where viscous forces dominate over inertial forces \cite{purcell}. Helical flagella allow noninvariant deformation under time-reversal symmetry, enabling displacement. This swimming mechanism is exemplified by the \textit{Escherichia Coli} bacterium, which naturally utilizes helical flagella for propulsion \cite{e_coli_2006}, or the \textit{MO-1} bacterium \cite{lefevre_M01, Zhang_M01}. These swimmers are also widely used in robotics because they can be easily guided by an external magnetic field \cite{Khalil2016}. The performance of microswimmers with single or dual helical flagella has been thoroughly studied, especially with ellipsoidal cell bodies, by optimizing certain parameters while keeping others fixed \cite{phan-thien_tran-cong_ramia_1987, higdon_helical_waves, two_flagella_shum}, using genetic algorithms \cite{kenta_optisperm}, or even experimentally \cite{mult_flagella_regnier}. With the development of Machine Learning (ML), new methods applied to these problems could yield interesting results, similar to those used in aerodynamics \cite{ml_aero_shapeopti}.\\

\noindent In this paper, we explore the optimization of microswimmer shapes with one and two helical flagella, aiming to find shapes that enable optimal swimming along a straight line. This involves maximizing speed while accounting for the power expended during movement. The challenge here is to optimize both the shape of the head and the flagella. This requires coupling various methods to solve the fluid-structure interaction and the optimization problem. The dynamics of the swimmer are obtained using the Boundary Element Method (BEM), commonly applied in such applications \cite{Pozrikidis_1992, WALKER2019311, two_flagella_shum, phan-thien_tran-cong_ramia_1987}. The BEM is implemented using the \texttt{MATLAB} library \texttt{Gypsilab} \cite{alouges2018bem}. Several shape optimization methods are available, including parametric, geometric, and topological optimization \cite{allaire2007conception, ALLAIRE20211}. We focus on parametric shape optimization due to the complexity of the model, particularly the challenges associated with computing the shape gradient \cite{moreau2022shapes}, especially when constraints are involved. The use of Bayesian Optimization (BO) methods, specifically the Scalable Constrained Bayesian Optimization (SCBO) method \cite{eriksson_scalable_2021}, represents a novel approach in microswimmer shape optimization. This method, implemented in \texttt{PYTHON} via the \texttt{BoTorch} library \cite{balandat2020botorch}, is particularly well-suited for high-dimensional problems with numerous constraints.\\

\noindent One of the major novelties of this paper is that the swimmer's head is not constrained to a particular geometric shape; instead, a generic shape can be obtained through optimization. Most previous studies have focused on ellipsoidal heads \cite{phan-thien_tran-cong_ramia_1987, gadelha_optimal_2013, kenta_optisperm, two_flagella_shum, chwang_helical_movement, higdon_helical_waves, keller_swimming_1976}. However, the shape of a microswimmer's head can be complex, such as the flattened water drop shape of human spermatozoa \cite{SMITH_2009}. We use the Free-Form-Deformation (FFD) approach \cite{ffd} to achieve a generic shape framework. FFD allows the deformation of a parameterized shape representation by manipulating control points, which is then coupled with BO methods, as demonstrated in \cite{ffd_bo_2023}.\\

\noindent The optimal shape results are validated and benchmarked against swimmers found in the literature \cite{phan-thien_tran-cong_ramia_1987, two_flagella_shum}. Simulations are performed for two optimization problems with different goals: maximize average speed or maximize average efficiency. We study both monoflagellated and biflagellated swimmers. New swimmer forms emerge from this paper, providing insights into microscale swimming dynamics. These forms are discussed and compared with their biological counterparts, showing a rich variety of swimmers that can be \textit{pushers} or \textit{pullers} \cite{Lauga_2009}.\\

\noindent The paper is organized as follows: \Cref{sec:1} presents the mathematical modeling of microswimmers, including their geometry and the swimming problem's numerical scheme. \Cref{sec:2} introduces the shape optimization problem with FFD and BO coupling. 
\Cref{sec:3} presents the numerical results, beginning with the optimization of a few parameters for swimmers found in the literature, and then extending the analysis to include the optimization of all parameters. \Cref{sec:4} and \Cref{sec:5} summarize and conclude the paper.

\section{Microswimmers modelling}\label{sec:1}
This section provides a description of the geometry of microswimmers, with a focus on the flagella, followed by an explanation of the system's dynamics and the numerical scheme used to solve them. The parameterization of the head is detailed in \Cref{sec:2}.
\subsection{Geometry of microswimmers}
\begin{figure}[htpb]
    \centering
    \begin{tikzpicture}
        \matrix[matrix of nodes,  column sep=2cm]{
        \includegraphics[width=6cm]{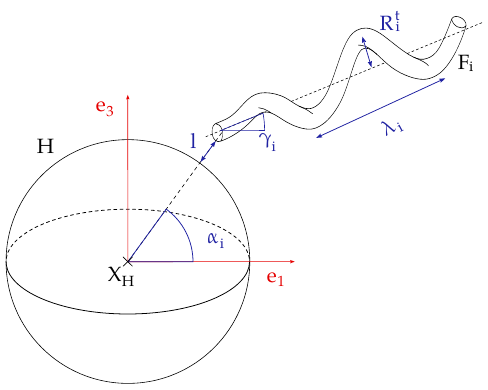}&  
       \includegraphics[width=6cm]{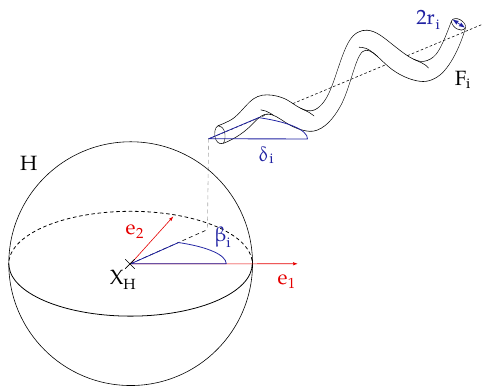}\\
        };
    \end{tikzpicture}
   \caption{Schematic of microswimmers represented by an arbitrary flagellum $F_i$ with parameters $(\alpha_i, \beta_i, \gamma_i, \delta_i, \lambda_i, R_i^t, r_i, l_i)$ in the frame $(e_1, e_3)$ (left) and in the frame $(e_1, e_2)$ (right).}
    \label{fig:geo_swimmer}
\end{figure}

The swimmer, denoted by $S$, is composed of a head $H$ and a set of $n_f$ of flagella, i.e.\ $S = H \bigcup_{i=1}^{n_F} F_i$. The head is assumed to be an arbitrary shape in $\mathbb{R}^3$. Each flagellum $F_i$ is a tube of radius $r$ with a centerline of total length $L$ described by, for $s \in [0,L]$
\begin{equation*}\label{eq:flag}
\begin{cases}
    x(s) &= s,\\
    y(s) &= R^t\left(1 - e^{-k_{E}^2 s^2}\right)\cos\left(\frac{2\pi s}{\lambda}\right),\\
    z(s) &= R^t\left(1 - e^{-k_{E}^2 s^2}\right)\sin\left(\frac{2\pi s}{\lambda}\right),
\end{cases}
\end{equation*}
where $R^t$ is the maximal amplitude of the helix, $\lambda$ is the wavelength, and $k_E$ is a shrinkage coefficient \cite{higdon_helical_waves}. The distance between each flagellum and the head is defined by $l$ along the normal to the head. We denote $(e_1, e_2, e_3)$ the body frame which moves with the swimmer, centered in the center mass of the head $x^H$. The junction $x^{F_i}$ of each flagellum $F_i$ is characterized by angles $\alpha_i$ in the frame $(e_1, e_3)$ and $\beta_i$ in the frame $(e_1, e_2)$. The orientation of each flagellum is determined by the angles $\gamma_i$ and $\beta_i$, representing rotations around $e_2$ and $e_3$ respectively. The geometric configuration of the swimmer is illustrated in \Cref{fig:geo_swimmer}, where the head is represented by a sphere and only one flagellum is depicted for simplicity. This modelisation was already employed in \cite{two_flagella_shum} for  biflagellated
swimmer, and in \cite{phan-thien_tran-cong_ramia_1987} for monoflagellated one.

\subsection{Dynamics of microswimmers}
Due to the low Reynolds number in microswimmer dynamics, the fluid behavior is governed by the Stokes equations. Let $u : \mathbb{R}^3 \to \mathbb{R}^3$ and $p : \mathbb{R}^3 \to \mathbb{R}$ represent the velocity and pressure fields of the fluid, respectively, satisfying the nonslip condition and the following equations:
\begin{equation}\label{eq:stokes}
    \left\{
    \begin{aligned}
    -\mu \Delta u + \nabla p &= 0 &\text{in } \mathbb{R}^3\setminus S, \\
    \nabla \cdot u &= 0 &\text{in } \mathbb{R}^3\setminus S, \\
    u &= U + \Omega \wedge (x-x^H) &\text{on } \partial H, \\
    u &= U + \Omega \wedge (x-x^H) + \omega e_1^{F_i} \wedge (x-x^{F_i}) &\text{on } \partial F_i,\\
    \|u\|, p &\to 0 &\text{as } x\to +\infty,
    \end{aligned}
    \right.
\end{equation}
for all $i \in \{1, \ldots, n_F\}$ and where $S$ denotes the microswimmer. Here, $\mu=1$ represents the viscosity of the fluid, and $\omega=-2\pi \text{rad.s}^{-1}$ denotes the angular velocity of the flagella around their directional axes $e_1^{F_i}$, mimicking the propulsion mechanism of bacteria that generate helical waves along their flagella. The motion of the rigid body is decomposed into a translational velocity $U \in \mathbb{R}^3$ and an angular velocity $\Omega \in \mathbb{R}^3$ induced by the rotation of the flagella of angular velocity $\omega$. The Dirichlet boundary conditions are imposed by $\omega e_1^{F_i} \wedge (x-x^{F_i})$, where $\omega$ is known, and $U + \Omega \wedge (x-x^H)$, where $U$ and $\Omega$ are the variables of interest that we aim to determine.

\subsection{Numerical scheme}
Stokes' equations, being linear with respect to velocity $u$ and $p$, admit the existence of a Green's tensor kernel $G_{ij}$ for $(i,j)\in\{1,2,3\}^2$. For $(x,y)\in\left(\mathbb{R}^3\right)^2$, the Green's function is given by:
\begin{equation*}
    G_{ij}(x,y)=\frac{1}{8\pi\mu}\left(\frac{\delta_{ij}}{\|x-y\|_2}+\frac{(x_i-y_i)(x_j-y_j)}{\|x-y\|_2^3}\right),
\end{equation*}
where $\delta_{ij}$ is the Kronecker delta. The convolution of the Green's kernel with the fluid surface tensions $f=\sigma(u,p)n$ yields an integral representation of the velocity field $u$ on $\partial S$. Thus, for any $y\in \partial S$ \cite{Pozrikidis_1992, happel_brenner_1983}, we have:
\begin{equation}\label{eq:ui_int}
    u_i(y) = -\int_{\partial S} \sum_{j=1}^3 G_{ij}(x,y)f_j(x)dx + \frac{1}{8\pi}\int_{\partial S}\sum_{j,k=1}^3 u_j(x)T_{jik}(x,y)n_k(x)dx,
\end{equation}
where $T_{jik}(x,y)=-6\frac{(x_i-y_i)(x_j-y_j)(x_k-y_k)}{\|x-y\|_2^5}$ for $(i,j,k)\in\{1,2,3\}^3$. In the Stokes equation \eqref{eq:stokes}, the immersed object describes the motion of a rigid object with velocity $u(x)=U+\Omega\wedge (x-x^H)$, leading to:
\begin{equation}\label{eq:sec_term_uint}
    \frac{1}{8\pi}\int_{\partial F}\sum_{j,k=1}^3 u_j(x)T_{jik}(x,y)n_k(x)dx=0.
\end{equation}
Substituting \eqref{eq:sec_term_uint} into \eqref{eq:ui_int}, we obtain:
\begin{equation}\label{eq:u_int}
    u(y)=-\int_{\partial S}G(x,y)f(x)dx.
\end{equation}
Equation \eqref{eq:u_int} provides the velocity expression on the solid boundary once the stress tensor is known. It can be inverted to obtain the stress tensor on the solid boundary. If the solid is subjected to external forces or torques, $F_{\text{ext}}$ and $T_{\text{ext}}$ respectively, \eqref{eq:u_int} must be complemented by the following equations:
\begin{equation}\label{eq:FT_equations}
\begin{aligned}
    \int_{\partial S}f(x)dx &= F_{\text{ext}},\\
    \int_{\partial S}f(x)\wedge (x-x^H)dx &= T_{\text{ext}}.
\end{aligned}
\end{equation}
If the external forces and torques are zero, the equations \eqref{eq:FT_equations} are called \textit{self-propulsion constraints}. Using \eqref{eq:stokes}, \eqref{eq:u_int}, and \eqref{eq:FT_equations} for self-propulsion, we obtain a system for determining the translational and rotational speeds $(U,\Omega)$, given by:
\begin{equation}\label{eq:stokes_boundary}
\begin{aligned}
     U-(x-x^H)\wedge\Omega + \int_{\partial S}G(x,y)f(y)dy &= 0, & \forall x &\in \partial H,\\
     U-(x-x^H)\wedge\Omega + \int_{\partial S}G(x,y)f(y)dy &= (x-x^{F^i})\wedge \omega e_1^{F^i}, & \forall x &\in \partial F^{i},\\
     \int_{\partial S}f(y)dy &=0,\\
     \int_{\partial S}f(y)\wedge (y-x^H)dy &= 0.
\end{aligned}
\end{equation}
Expressing the Stokes problem \eqref{eq:stokes}  as a problem on $S$ boundaries \eqref{eq:stokes_boundary}, we can use the BEM. In our case, the computational domain is simply the mesh that discretizes the swimmer, i.e., the mesh of the head and flagella. Let $N$ be the number of nodes in the mesh and $V_N$ the conformal finite element space we're working in. Here, we will consider $\mathbb{P}^1$ finite elements. Expressing $f$ in the $V_N$ basis, we have for all $y \in \partial S$ and for all $j\in\{1,2,3\}$:
\begin{equation*}
    f_j(y)=\sum_{l=1}^N f_j^l \phi^l(y),
\end{equation*}
where $\phi=\{\phi^1,\ldots, \phi^N\}$ is the basis of $V_N$ on $\partial S$ and $\begin{pmatrix}f_j^1 &\ldots &f_j^N\end{pmatrix}^T$ is the coordinate vector of $f_j$ in $V_N$. By denoting $f^l=\begin{pmatrix}
    f_1^l & f_2^l & f_3^l
\end{pmatrix}^T$ and multiplying by $\phi^k$, an element of $V_N$, followed by integration over $\partial S$ of the first two equations in \eqref{eq:stokes_boundary}, for $k \in \{1, \ldots, N\}$, we obtain following the system:
\begin{equation}\label{eq:bem_stokes}
\begin{aligned}
     \int_{\partial H}\left(U-(x-x^H)\wedge\Omega\right) \phi^k(x)dx &+ \int_{\partial H}\left(\int_{\partial S} G(x,y)\sum_{l=1}^N f^l \phi^l(y)dy\right)\phi^k(x)dx \\
      &= 0, \\
     \int_{\partial F^i}\left(U-(x-x^H)\wedge\Omega\right) \phi^k(x)dx &+ 
     \int_{\partial F^i}\left(\int_{\partial S}G(x,y)\sum_{l=1}^N f^l \phi^l(y)dy\right)\phi^k(x)dx \\&= \int_{\partial F^i}(x-x^{F^i})\wedge \omega e_1^{F^i}\phi^k(x)dx, \\
     \int_{\partial S}\sum_{l=1}^N f^l \phi^l(y)dy &=0,\\
     \int_{\partial S}(y-x^H)\wedge \sum_{l=1}^N f^l \phi^l(y)dy &= 0,
\end{aligned}
\end{equation}
for all $i \in \{1, \ldots, n_F\}$. Since $\displaystyle{S = H \bigcup_{i=1}^{n_F} F_i}$, we can separate the integrals so that \eqref{eq:bem_stokes} reads as a matrix system:

\begin{equation*}
\begin{bmatrix}
    G & J^T & K^T \\
    J & 0 & 0\\
    K & 0 & 0
\end{bmatrix}\begin{bmatrix}
    f\\
    U\\
    \Omega 
\end{bmatrix}=\begin{bmatrix}
    I(\omega)\\ 0 \\0
\end{bmatrix}.
\end{equation*}
Submatrix $G$ is composed of $(1+n_F)^2$ submatrices $G^{\{A,B\}}$ for $(A,B)\in \{\partial H, \partial F^1, \ldots, \partial F^{n_F}\}$. Each $G^{\{A,B\}}$ matrix can be subdivided into $N_A \times N_B$ submatrices of size $3\times 3$ defined by:
\begin{equation*}
    G^{\{A,B\}}_{ij,lk} = \int_{A}\int_{B}G_{ij}(x,y)\phi^l(x)\phi^k(y)dxdy,
\end{equation*}
for $(i,j)\in\{1,2,3\}^2$, $l\in\{1,\ldots, N_B\}$, and $k\in\{1,\ldots, N_A\}$. The submatrix $J$ is composed of $(1+n_F)$ submatrices $J^A$ for $A\in\{\partial H,\partial F^1, \ldots,\partial F^{n_F}\}$. Each matrix $J^A$ is a diagonal block matrix where each block is of dimension $1\times N_A$ and the $j$-th is denoted $J^A_{il,j}$, defined by:
\begin{equation*}
    J^A_{il,j} = \int_{A}e_i \phi^l(x)dx,
\end{equation*}
for $(i,j)\in\{1,2,3\}^2$ and $l\in\{1,\ldots, N_A\}$. Similarly, the submatrix $K$ is composed of $(1+n_F)$ submatrices $K^A$ for $A\in\{\partial H,\partial F^1, \ldots,\partial F^{n_F}\}$. Each matrix $K^A$ has size $3\times 3N_A$, and its components $K_{ij,l}^A$ are defined by:
\begin{equation*}
    K^A_{ij,l} = \int_{A}\left[(x-x^H)\wedge e_i\right]_j \phi^l(x)dx,
\end{equation*}
for $(i,j)\in\{1,2,3\}^2$ and $l\in\{1,\ldots, N_A\}$. The vector $I(\omega)$ is composed of $(1+n_F)$ subvectors $I(\omega)^{\partial H}=0$ and $I(\omega)^A$ for $A \in \{\partial F^1, \ldots, \partial F^{n_F}\}$. Each vector $I(\omega)^A$ is subdivided into $N_A$ subvectors of size 3, defined by:
\begin{equation*}
    I(\omega)_{l}^A=\int_{A}\left[(x-x^H)\wedge \omega e_1^A\right]\phi^l(x)dx,
\end{equation*}
for $l\in\{1,\ldots, N_A\}$. The implementation is done using the open-source BEM code \texttt{Gypsilab}\footnote{\url{https://github.com/matthieuaussal/gypsilab}} in Matlab and the validation cases are presented in \Cref{app:valBEM}.

\section{Shape Optimization Problem}\label{sec:2}
We introduce $S^0$ as the reference monoflagellated swimmer, whose characteristics are summarized in \Cref{table:S0}, and whose flagellar parameters correspond to those in \cite{two_flagella_shum}, for future comparison. The resulting average velocities and power dissipation are tabulated in \Cref{table:S0_res}. 
\begin{table}[htpb]
\caption{\label{table:S0} Parameter values for the cell body and flagellum of the reference swimmer $S^0$.}
\begin{ruledtabular}
    \begin{tabular}{l c c}
         Parameter & Value (Dimensionlesss) & Value (Dimensional) \\
         \colrule
         $H^0$ & $B_{\|\cdot\|_2}(0,1)$ & $B_{\|\cdot\|_2}(0,0.74)$ \\ 
         $L^0$ & $3$ & $2.2$ µm \\
         $r^0$ & $0.067$ & $0.05$ µm \\  
         $\lambda^0$ & $1$ & $0.74$ µm \\ 
         $R^{t0}$ & $0.2$  & $0.15$ µm\\ 
         $l^0$ & $2r^0$ & $0.1$ µm \\
         $k_E^0$ & $0.333\times 2\pi/\lambda^0$ & $2.8$ $\text{µm}^{-1}$\\
              $\alpha^0,\gamma^0, \beta^0, \delta^0$ & $0$ & $0$ rad  \\ 
    \end{tabular}
\end{ruledtabular}   
\end{table}
\begin{figure}[htpb]
    \centering  \includegraphics[width=7cm]{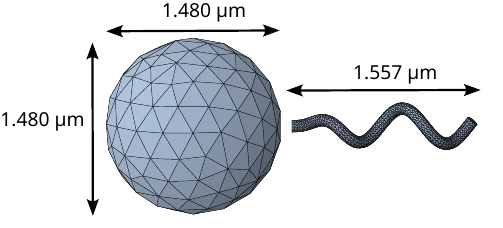}
    \caption{Discretization of the reference swimmer $S^0$ used to obtain values in \Cref{table:S0_res}.}
    \label{fig:s0}
\end{figure}

\subsection{Admissible shape}
In this section, we define the set of shapes from which potential swimmers are selected.
As outlined earlier, each swimmer comprises two components: the head and its flagella.
These components belong to respective sets of \textit{admissible shapes}, denoted as $\mathcal{S}_{ad} = \mathcal{H}^{\varepsilon}_{ad} \times \mathcal{F}^{\mathcal{P}}_{ad}$, where $\mathcal{F}^{\mathcal{P}}_{ad}$ (for flagella) and $\mathcal{H}^{\varepsilon}_{ad}$ (for the head) are defined subsequently.
In what follows, we denote by $|A|$ the volume of the $A$ shape.

\subsubsection{Flagella admissible shape}
We define the set of flagella admissible shapes $\mathcal{F}^{\mathcal{P}}_{ad}$  by
\begin{equation}\label{eq:flag_ad}
    \mathcal{F}^{\mathcal{P}}_{ad} = \left\{~F \subset \mathbb{R}^3\mid  (\lambda, R^t, \alpha, \gamma, \beta, \delta)\in \mathcal{P} \quad \text{and} \quad |F|=|F^0|~\right\},
\end{equation}
with $\mathcal{P}$ a compact set associated with the problem under study and $F^0$ the referent flagella whose characteristics are shown in \Cref{table:S0}. To keep flagellum volumes constant, we will consider flagella of the same radius $r=r^0$ and the same total length $L=L^0$.

\subsubsection{Head admissible shape}
In this paper, we employ FFD framework \cite{ffd} to parametrize the shape of the swimmer's head. The FFD approach describes an arbitrary basic shape using a discrete set of control points. Let's briefly outline the main framework below.\\

\begin{figure}[htpb]
\centering
\def\svgwidth{10cm}
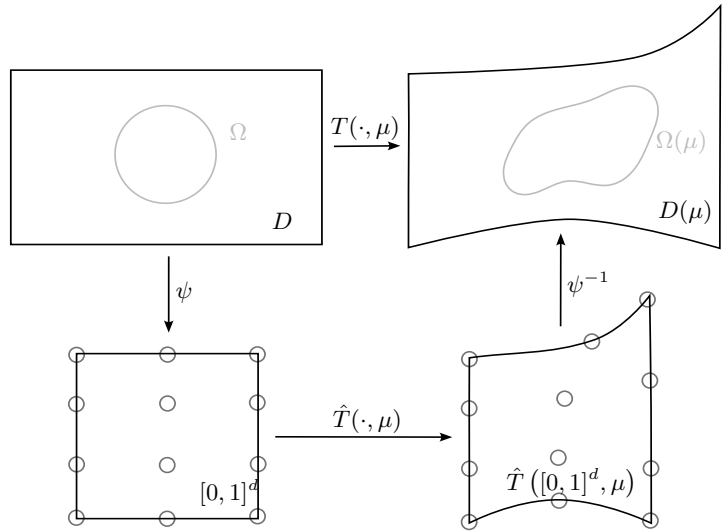
\caption{Illustration of the FFD method deforming a domain $D$ containing a $\Theta$ shape via a $\mu$ vector using the $T$ application described by \eqref{eq:ffd}. }
\label{fig:ffd}
\end{figure}

\noindent The shape we are working on, denoted by $\Theta$, is embedded in a slightly larger geometric domain denoted by $D$.
The key point is to operate within the unit cube of $\mathbb{R}^3$ by using a diffeomorphism with $D$. In this unit cube, $M$ points, called the control points, are defined to track the geometrical deformation of the shape. Moving these control points induces a deformation of the cube consistent with their displacement.
The resulting deformed shape is obtained by applying the inverse diffeomorphism to the ''deformed unit cube''(see \Cref{fig:ffd}). 
This setup defines a map, denoted by $T$, which deforms the reference shape by shifting the control points by an amount $\mu \in \mathbb{R}^{3M}$,
\begin{equation}\label{eq:ffd}
    \begin{array}{rcl}
            T :D \times \mathbb{R}^{3 \times M} &\to& D(\mu)\\
            (\Theta, \mu) &\mapsto & \left(\psi^{-1}\circ \hat{T} \circ \psi\right)(\Theta,\mu),
            \end{array}
\end{equation}
as illustrated in \Cref{fig:ffd}. We define the set from which we can select the potential heads as follows:
\begin{equation}\label{eq:head_ad}
   \mathcal{H}_{ad}^{\varepsilon} = \left\{ H = T(H^0, \mu)\ \middle\vert \begin{array}{l}
    \mu \in \mathcal{V} \quad \text{and} \quad \left||H|-|H^0|\right|\leq \varepsilon
  \end{array}\right\}.
\end{equation}
Here, the volume of the sphere $H^0$ is preserved with a precision of $\varepsilon$ and $\mathcal{V}$ defines a sufficiently small closed set around the control points to prevent mesh collapse. For $\mathcal{V}$, the displacement of each control point is constrained to prevent excessive overlap or crossing between points, while only the boundary points of the domain are allowed to move, keeping the internal points fixed. In \Cref{app:ffdillu}, several shapes are illustrated by applying the FFD map to the reference head \( H^0 \) with different control points. These shapes, which lie in \( \mathcal{H}_{ad}^{\varepsilon} \), are shown for various numbers of control points. By adjusting the control points, we can obtain a wide variety of shapes, including nonconvex ones.

\subsection{Statement of the optimization problems}
In this part, we state the optimization problem by defining the cost functions of interest. First, we define the power dissipated and the average power dissipated as follows \cite{two_flagella_shum}:
\begin{equation}\label{eq:mean_P}
    P = \sum_{i=1}^{n_F} -\omega e_1^{F_i} \cdot \int_{\partial F_i} \sigma(u,p)n \wedge (x-x^{F_i}) \quad \text{and} \quad \bar{P} = \frac{1}{T}\int_{0}^T P(t)dt,
\end{equation}
where $T$ is the time for a single stroke, $n$ is the outward normal to $\partial F$,  $\sigma$ denotes the stress tensor, given by $\sigma(u,p)=-pI + 2\mu e(u)$, and $e(u)$ represents the deformation tensor, defined as $e(u)=\frac{1}{2}(\nabla u + \nabla u^T)$. 
Then, the mean translational and angular velocities are given by:
\begin{equation}\label{eq:mean_vel}
    \bar{U} = \frac{1}{T}\int_{0}^T U(t)dt \quad \text{and} \quad \bar{\Omega} = \frac{1}{T}\int_{0}^T \Omega(t)dt.
\end{equation}
\begin{table}[htpb]
\caption{\label{table:S0_res} Mean velocities and power dissipated obtained from \eqref{eq:stokes} by taking the swimmer $S^0$ from \Cref{table:S0}.}
\begin{ruledtabular}
  \begin{tabular}{l c c}
    Parameter & Value (Dimensionless) & Value (Dimensional)\\
\colrule
      $\bar{U}^0_1$ & $-0.0306$ & $-0.0226$ $\text{µm.s}^{-1}$ \\ 
     $\bar{\Omega}^0_1$& $0.1470$ &  $0.1470$ $\text{rad.s}^{-1}$\\
     $\bar{P}^0$&$-24.7512$ & $-18.276$ $\text{µW}$\\
  \end{tabular}
\end{ruledtabular}
\end{table}
Finally, we define the cost functions, $J_1$ and $J_2$, corresponding respectively to the mean speed in the direction of $e_1$ normalized by the one of the reference swimmer $S^0$ and the mean efficiency, as:
\begin{equation}\label{eq:cost_func}
    J_1(S) = -\frac{\bar{U}_1(S)}{\bar{U}^0_1} \quad \text{and} \quad J_2(S) = -\frac{\bar{U}_1(S)}{\bar{U}^0_1} \times \frac{\bar{P}^0}{\bar{P}(S)}.
\end{equation}
\noindent Hence, the optimization problem, for $i \in \{1,2\}$, is formulated as:
\begin{equation}\label{eq:opti_pb}
    \inf_{S \in \mathcal{S}_{ad}} J_i(S).
\end{equation}
Cost functions are adapted to work with minimization problems, as it is more conventional to address minimization rather than maximization in optimization. In the following, to avoid notational complexity, we omit the dependence on the swimmer's shape in some functional.

\subsection{Bayesian optimization Applied to shape optimization}

\begin{figure}[h]
    \centering
    \includegraphics[width=12cm]{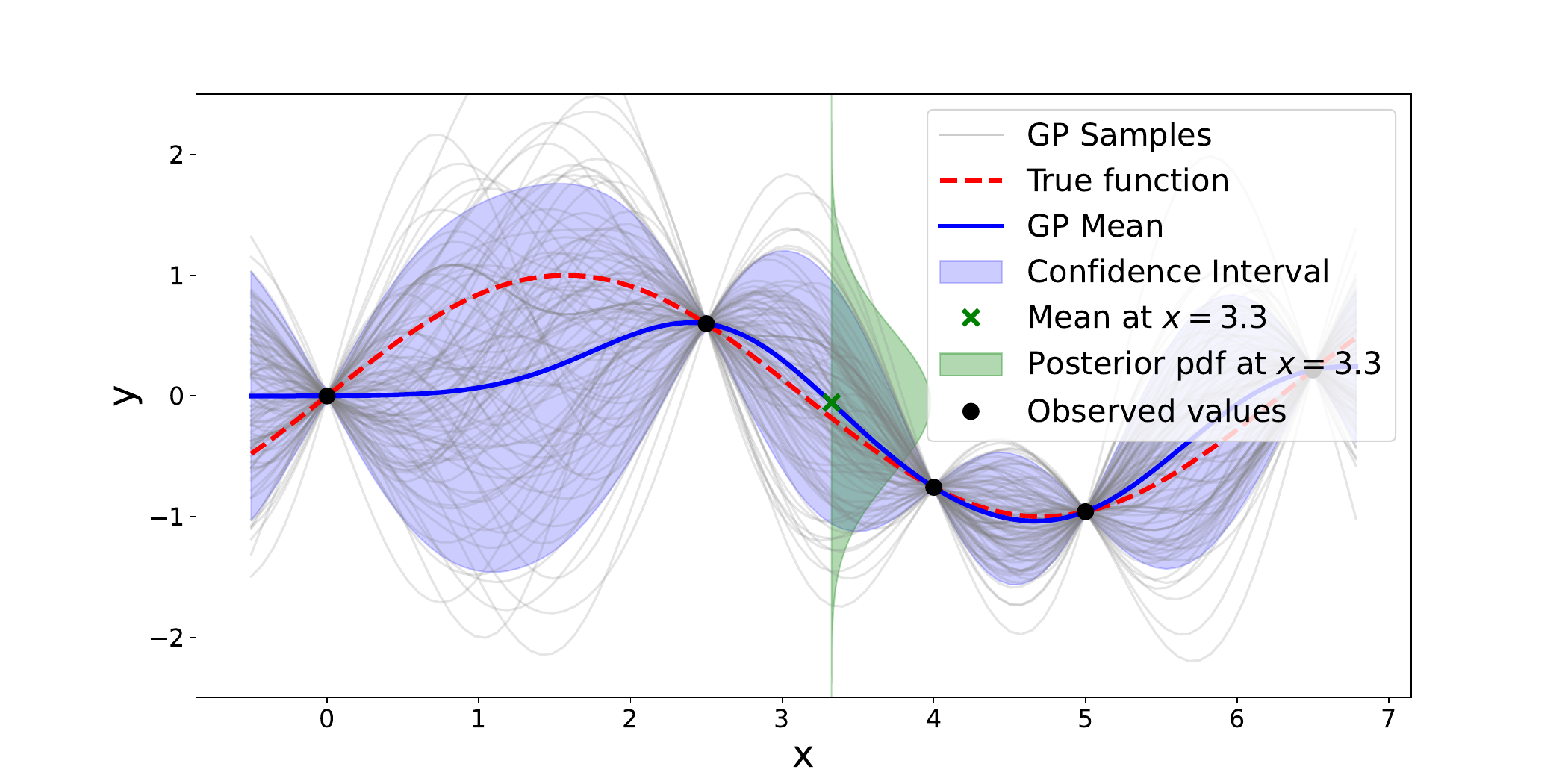}
    \caption{Gaussian process on a toy function. The true function (red line) is approximated by the GP mean (blue line) according the observations (black point). The Gaussian posterior distribution for $x=3.3$ is depicted (green). The associated point-wise confidence intervals at level 95\% are displayed (blue). A set of $100$ GP samples is displayed (gray).}
    \label{fig:gp}
\end{figure}

To solve these microswimmer shape optimization problems, we apply Bayesian optimization methods \cite{garnett2023bayesian}, specifically the Scalable Constrained Bayesian Optimization (SCBO) method \cite{eriksson_scalable_2021}. First of all, Gaussian processes (GPs) are probabilistic modeling techniques used to approximate functions when direct evaluation is computationally expensive \cite{gramacy2020surrogates, Rasmussen2006}. Instead of evaluating the expensive function $J$ for many swimmer shapes, GPs allow us to model $J$ based on limited data and make predictions about its behavior with uncertainty estimates. \Cref{fig:gp} illustrates how GPs work on a toy function, showing observed values and the associated GP model.

\begin{figure}[htpb]
\centering
\def\svgwidth{14cm}
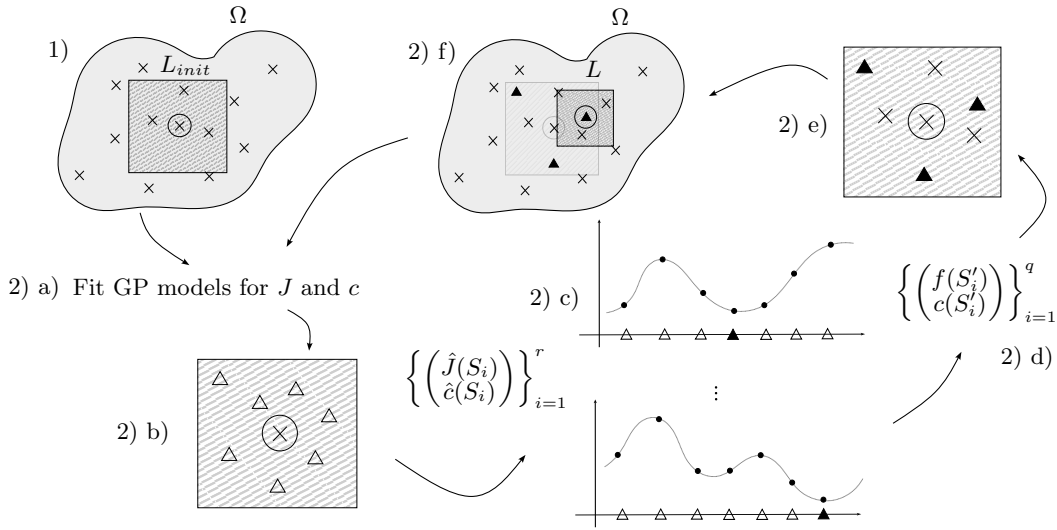
\caption{SCBO representation: (1) Initial shapes are taken from $\Omega$ (represented by crosses), and the trust zone is defined (a grid with the central point surrounded and the edge of the zone defined by a square). (2a) The Gaussian process models associated with the objective $J$ and the constraints $c$ are constructed. (2b) $r$ random shapes (represented by triangles) are taken from the trust region. (2c) A batch of $q$ realizations $\{\begin{pmatrix}\hat{J}(S_i) & \hat{c}(S_i)\end{pmatrix}^T\}_{i=1}^r$ is calculated. For each realization, the $x$ axis corresponds to the $S$ shapes in $\Omega$ and is discretized by the $\{S_i\}_{i=1}^r$ (the triangles) chosen previously. The $y$ axis shows the realizations. For simplicity, a single graph is drawn for the objective and the constraints. For each realization, we keep the best shape (completely filled triangle). (2d) For the $q$ best shapes $\{S_i'\}_{i=1}^q$, we calculate the real value of the objective and the constraints. (2e) These $q$ new shapes are added to the previous observations. (2f) The trust region is readjusted by modifying its center and length.}
\label{fig:scbo}
\end{figure}

\noindent Solving the problem \eqref{eq:opti_pb} requires rewriting the constraint spaces $\mathcal{S}_{ad}$, defined by \eqref{eq:flag_ad} and \eqref{eq:head_ad}, in terms of inequality constraints such as:
\begin{equation*}
\mathcal{S}_{ad}=\left\{S=(H,F)\subset \mathbb{R}^3 \times \mathbb{R}^3 \mid c(S)=(c^1(S), \ldots, c^m(S))\leq 0\right\}.
\end{equation*}
The main method is presented above with \Cref{fig:scbo}, which depicts the different steps.
\begin{enumerate}
\item We start by taking a number of initial shapes and defining the trust region using maximum utility. The trust region approach allows the selection of samples locally around the best current point by defining a hypercube. Given $n$ shapes $\{S_i\}_{i=1}^{n}$, we evaluate the cost function $J$ and the constraint function $c$. Then, we define the set of feasible constraints for these shapes as:
\begin{equation*}
F_c=\left\{S_i \mid c(S_i) \leq 0\right\}.
\end{equation*}
To choose the center $S^*$ of the hypercube, two cases are considered depending on the cardinality of $F_c$:
\begin{equation*}
S^* = \begin{cases}
\displaystyle{\argmin_{S \in F_c} J(S)} \quad &\text{if } F_c \neq \emptyset, \\
\displaystyle{\argmin_{S \in \{S_i\}_{i=1}^n} \sum_{j=1}^m \max\left\{c^j(S), 0\right\}} \quad &\text{if } F_c = \emptyset,
\end{cases}
\end{equation*}
i.e., either $F_c$ is nonempty, and we take the point minimizing the cost function, or $F_c$ is empty, and we take the point minimizing the maximum constraint violation. The length of the hypercube is initialized by $L_{\text{init}}$.

\item Until the number of remaining possible evaluations reaches zero or $L<L_{\text{min}}$:
\begin{enumerate}
\item Gaussian process models associated with the cost function and the constraint functions are constructed from the observations, i.e., shapes and evaluated functions for these shapes.
\item Randomly choose $r$ shapes inside the hypercube.
\item We consider a number of $q$ batches. For each batch, we obtain a sample from the posterior distribution of cost and constraints models (as in gray in \Cref{fig:gp}) at the level of discretization obtained by the previous $r$ shapes. The best of the $r$ shapes is then selected for each batch.
\item For the $q$ best shapes obtained in the previous step, the cost function and constraints are evaluated.
\item These new shapes are added to the previous observations.
\item The center of the trust region is updated based on the improvements achieved by the new candidates. Specifically, the new center is set as the current best point.  The length of the trust region is adjusted based on the number of successes or failures. Let $\tau_s$ and $\tau_f$ be the fixed success and failure rates, respectively, and $n_s$ and $n_f$ the number of successes and failures. Success is achieved when a better shape than the center of the hypercube is obtained. We say that $\tilde{S}$ is a better shape than $S$ when $S$ does not satisfy the constraints, and $\tilde{S}$ has a smaller constraint violation, or when $S$ satisfies the constraints, and $\tilde{S}$ also satisfies them while having lower objective value. The length of the hypercube is then adapted as follows:
\begin{equation*}
\text{if } n_s = \tau_s \text{ then } \begin{cases}
L = \min\{2L, L_{\text{max}}\},\\
n_s = 0,
\end{cases}
\text{ and }
\text{if } n_f = \tau_f \text{ then } \begin{cases}
L = L/2,\\
n_f = 0,
\end{cases}
\end{equation*}
with $L_{\text{max}}$ and $L_{\text{min}}$ denoting the maximum and minimum length of the hypercube.
\end{enumerate}
\item The algorithm returns the best admissible point found.
\end{enumerate}

\noindent A disturbance probability, as introduced in \cite{eriksson_scalable_2021}, has been incorporated to mitigate the edge effect, which refers to the phenomenon where points in high-dimensional spaces tend to lie near the boundaries of the search space. This effect can trap the algorithm at the boundaries of the hypercube. The implementation is done using the open-source SCBO code \texttt{BoTorch}\footnote{\url{https://github.com/pytorch/botorch}} in \texttt{PYTHON}.

\section{Numerical Results}\label{sec:3}
This section is divided into two parts. The first part focuses on presenting the results achieved using our method, which enhances the optimal results found in the literature \cite{phan-thien_tran-cong_ramia_1987,two_flagella_shum} constrained by a more restrictive geometrical framework. The second part addresses the outcomes obtained by optimizing all aspects for both the monoflagellated and biflagellated swimmers. Because of the nature of the problem, it is natural to consider two planes of symmetry for the head. This is described in \Cref{app:Numimpl}, along with the details of the numerical implementation. Detailed figures of the swimmers are given in \Cref{app:littswimmers} and \Cref{app:bestswimmer}. Bayesian optimization convergence graphs associated with optimal swimmers are presented in \Cref{app:cvgraph}.

\subsection{Outclassing competing models}
In this subsection, we concentrate on the top-performing swimmers identified in the case of monoflagellates from \cite{phan-thien_tran-cong_ramia_1987} and in the case of biflagellates from \cite{two_flagella_shum}. Using our method, we propose new swimmer designs capable of outperforming the previous models by relaxing the geometric assumptions about the swimmer's head. \\

\noindent Our results are presented in \Cref{table:outclassing_models}. The first row shows meshes of the swimmer we aim to outperform. In the first case, a monoflagellated swimmer from \cite{phan-thien_tran-cong_ramia_1987} is depicted whose characteristics are described in \Cref{table:87_param}, where the head is modeled as an ellipse defined by the equation
\begin{equation}\label{eq:ellipsoidal_head}
H=\left\{(x,y,z) \in \mathbb{R}^3 \mid \frac{x^2}{(R_1)^2}+\frac{y^2}{(R_2)^2}+\frac{z^2}{(R_3)^2}=1\right\},
\end{equation}
whose parameters have been optimised to minimize inverse efficiency given by
\begin{equation}\label{eq:inv_eff}
    \eta_{0}^{-1}=\frac{\bar{P}}{6\pi\mu\bar{A}\bar{U}_1^2},
\end{equation}
where $\bar{A}$ is the volume average radius of the cell body. Beside, we illustrate the swimmer designed using our method, which optimizes the head shape while maintaining the same head's volume and same flagellum by using the same cost function, i.e., inverse efficiency, given by \eqref{eq:inv_eff}. The resulting form is shorter but wider. Non-intuitively, the forward fluid-exposed side is broader and tapers towards the point where the flagellum is attached. This is contrary to expectations of a more streamlined, car-like front. The second row presents the case of a biflagellated swimmer. The position of the second flagellum is determined by a $\pi$ rotation of the first flagellum around the $e_1$ propulsion axis. This assumption allows us to constrain the swimmer's movement to the $e_1$ axis. The left swimmer, with a fixed ellipsoidal head (whose equation is given by \eqref{eq:ellipsoidal_head} and characteristics are described in \Cref{table:geo_shum}), has its junction between the head and the flagella optimized according to the mean efficiency $J_2$; see \cite{two_flagella_shum}. The swimmer depicted on the top right represents the result by optimizing only the head's shape. Subsequently, we optimized both the shape of the head and the angle of the junction between $0$ and $\pi/2$ as in \cite{two_flagella_shum}. For the swimmer whose head alone has been optimized, the shape is flatter and more elongated than the original. Unlike the monoflagellated swimmer, a triangular shape tapering towards the front can be seen. Additionally, when the angle of the junction between the head and the flagella is also optimized, the head shape becomes more elongated, echoing optimal aerodynamic shapes.\\ 

\noindent To compare the various microswimmers, we present their swimming characteristics in \Cref{table:num_res_outclass}. For the monoflagellated swimmer, the optimized version exhibits greater speed but consumes more power compared to the benchmark. Additionally, its wider, flatter geometry results in a significant reduction in angular velocity. In contrast, the optimized biflagellate swimmers demonstrate a higher angular velocity and lower power consumption. They also achieve greater translational speed.
\begin{table}[htpb]
\caption{\label{table:outclassing_models}
Comparison between benchmark swimmers and optimized designs based on efficiency metrics. The first column includes benchmark swimmers obtained from \cite{phan-thien_tran-cong_ramia_1987} (first row) and \cite{two_flagella_shum} (second row). The second column showcases outperforming designs achieved through different optimization strategies: the first and second rows represent designs optimized only by head shape, while the third row includes designs with the addition of the $\alpha$ parameter for biflagellate swimmers. Efficiency metrics are evaluated using inverse efficiency \eqref{eq:inv_eff} for monoflagellated swimmers and mean efficiency \eqref{eq:cost_func} for biflagellate swimmers.
}
\begin{ruledtabular}
  \begin{tabular}{ccc}
     & Benchmark  swimmers & Outclassing swimmers \\ 
  \colrule
     \raisebox{8.5ex}{Monoflagellated swimmers} & \includegraphics[width=5cm]{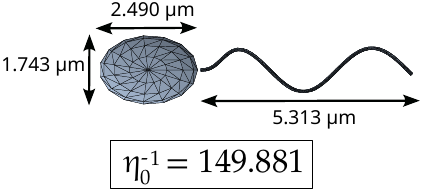}   &    \includegraphics[width=5cm]{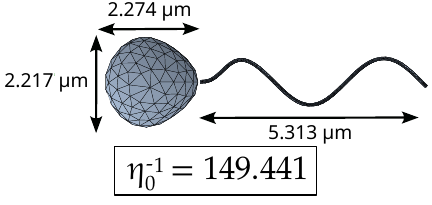}\\
  \colrule
     \multirow{1}*{\centering Biflagellate swimmers} &  \multirow{2}*[4.5em]{\includegraphics[width=5cm]{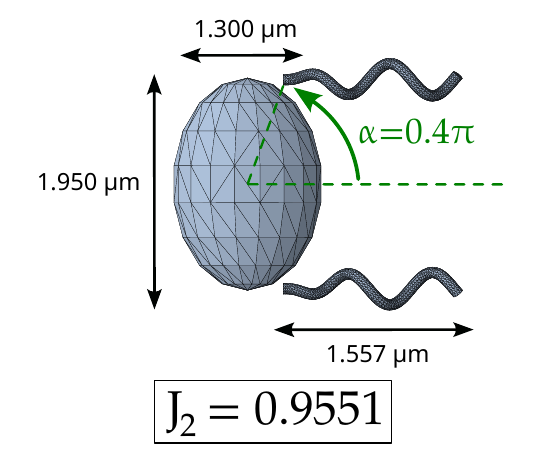}}&  \includegraphics[width=4.3cm]{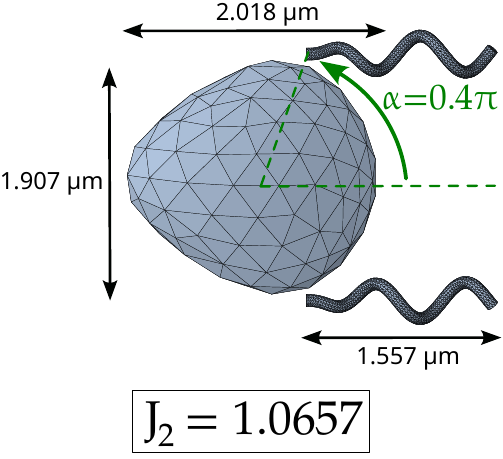} \\
     & &\includegraphics[width=5cm]{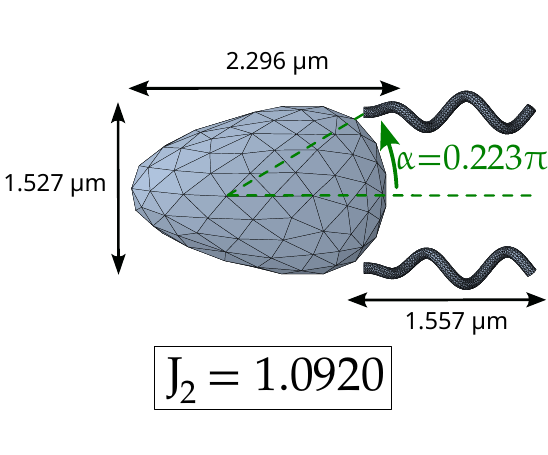}\\
  \end{tabular}
\end{ruledtabular}
\end{table}
\begin{table}[htpb]
\caption{\label{table:num_res_outclass} Swimming characteristics for benchmark swimmers from \cite{phan-thien_tran-cong_ramia_1987, two_flagella_shum} and optimized versions from \Cref{table:outclassing_models}. The swimming characteristics are normalized relative to the benchmark.}
\begin{ruledtabular}
  \begin{tabular}{c c c| c c c}
     & \includegraphics[width=2cm]{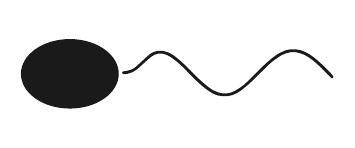} & \includegraphics[width=2cm]{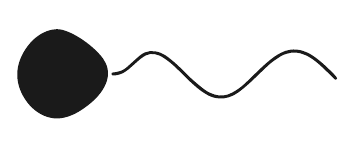} & \includegraphics[width=1cm]{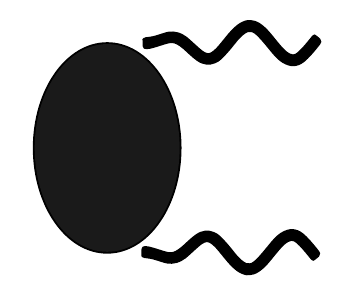}&  \includegraphics[width=1.1cm]{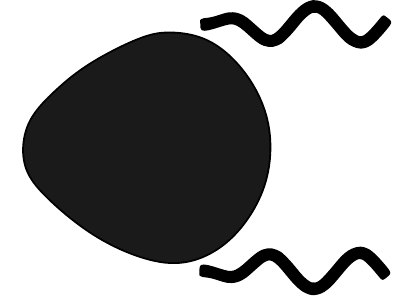}& \includegraphics[width=1.2cm]{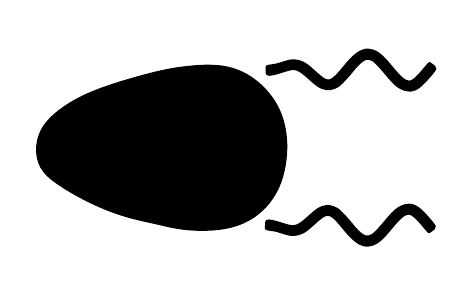}\\ 
     \colrule
     &&&&&\\[-1.5ex]
     $\bar{U}_1/\bar{U}^0_1$ & $1$ & $1.0618$ & $1$ & $1.1152$ & $1.1338$\\[0.5ex]
     $\bar{\Omega}_1/\bar{\Omega}^0_1$ &$1$ & $0.7983$ & $1$ & $1.0251$ & $1.3268$\\[0.5ex] 
     $\bar{P}/\bar{P}^0$  & $1$ & $1.1209$ &  $1$ & $0.9994$ & $0.9916$ \\ 
  \end{tabular}
\end{ruledtabular}
\end{table}

\subsection{Shape optimization of flagellated microswimmers}
In this section, we generalize the optimization process to account for all parameters describing the head and flagella of the microswimmer. We assume the swimmer's shape (head and the second flagellum in respect the first for biflagellate swimmers) must be symmetric with respect to $e_1$ (i.e., the desired direction of displacement). The flagella parameters are considered within the set $\mathcal{P}$, defining the admissible flagella $\mathcal{F}_{ad}^{\mathcal{P}}$ (see \eqref{eq:flag_ad}), prescribed as:
\begin{equation*}
    \mathcal{P} = [0.3, 4] \times [0.1, 1] \times \left[-\frac{\pi}{2}, \frac{\pi}{2}\right] \times \left[-\frac{\pi}{2}, \frac{\pi}{2}\right] \times \left[-\frac{\pi}{2}, \frac{\pi}{2}\right] \times \left[-\frac{\pi}{2}, \frac{\pi}{2}\right].
\end{equation*}
The existence and convergence of our approach towards a globally optimal shape are not confirmed. Therefore, multiple simulations of the same problem were conducted, consistently yielding similar results. The best results are summarized in \Cref{table:best_models} and \Cref{table:num_res}. Convergence graphs, presented in \Cref{app:cvgraph}, and detailed figures of the swimmers, provided in \Cref{app:bestswimmer}, offer additional insights.

\subsubsection{Monoflagellated swimmer}

\noindent {\bf Optimal shape description :} Focusing on the translational velocity, $J_1$ in \eqref{eq:cost_func}, \Cref{table:best_models} and \Cref{table:num_res} reveal that a water-drop-shaped head with a large-amplitude flagellum results in a swimmer about 3.6 times faster than the reference swimmer. When power dissipation is considered, $J_2$ in \eqref{eq:cost_func}, the results show a swimmer 1.5 times more efficient than the reference swimmer, with a slightly elongated head and a flagellum of smaller amplitude and longer wavelength. The head shapes are quite similar in both cases, so we will refer to these swimmers as \textit{water-drop swimmers}.
\\

\noindent {\bf Swimming trajectory :} As can be seen in \Cref{fig:best_traj}, the velocity-optimized swimmer has a helical trajectory with larger steps but also with a larger radius than the optimal swimmer that takes power dissipation into account. This is due to the large amplitude of the flagellum, which almost reaches the upper limit. Conversely, when power dissipation is considered, the resulting helical trajectory has a smaller radius. Thus, accounting for power dissipation significantly impacts the resulting optimal swimmer displacement. To meet the constraint specified in \eqref{eq:e1constaint} for average displacement along the $e_1$ axis, the algorithm identifies optimal swimmers with a flagellum oriented along the $e_1$ axis. However, a single flagellum induces significant displacement in the $e_2$ and $e_3$ directions during its stroke.\\

\subsubsection{Biflagellated swimmer}

\noindent {\bf Optimal shape description :} Focusing on the translational velocity, $J_1$ in \eqref{eq:cost_func}, we obtain a rather unusual swimmer with a much more compact head than those that can be obtained in other types of problem and flagella of very large amplitudes (here again the upper limit is almost reached). Because of its similarity in shape, we'll call this swimmer \textit{Bullhead swimmer}. This swimmer is almost 11 times faster than the reference swimmer $S^0$. However, as soon as we take into account the power expended, $J_2$ in \eqref{eq:cost_func}, we obtain a more "conventional" shape. It has an elongated head with flagella of low amplitude. The swimmer is slightly more than 1.2 times more efficient than the reference swimmer.\\

\begin{figure}[htpb]
    \centering
    \begin{tikzpicture}
        \matrix[matrix of nodes,  column sep=1cm]{
          \includegraphics[width=7.5cm]{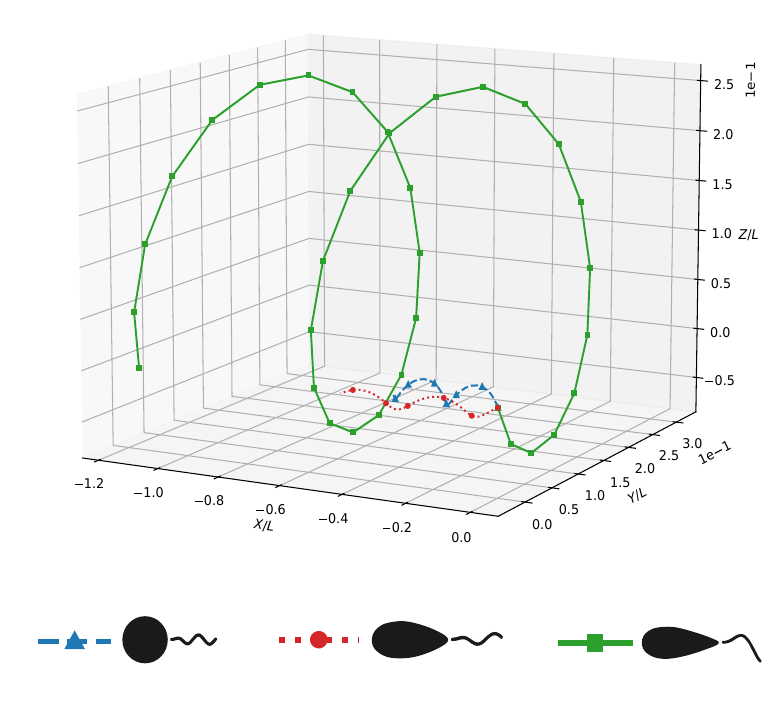} &
          \includegraphics[width=7.5cm]{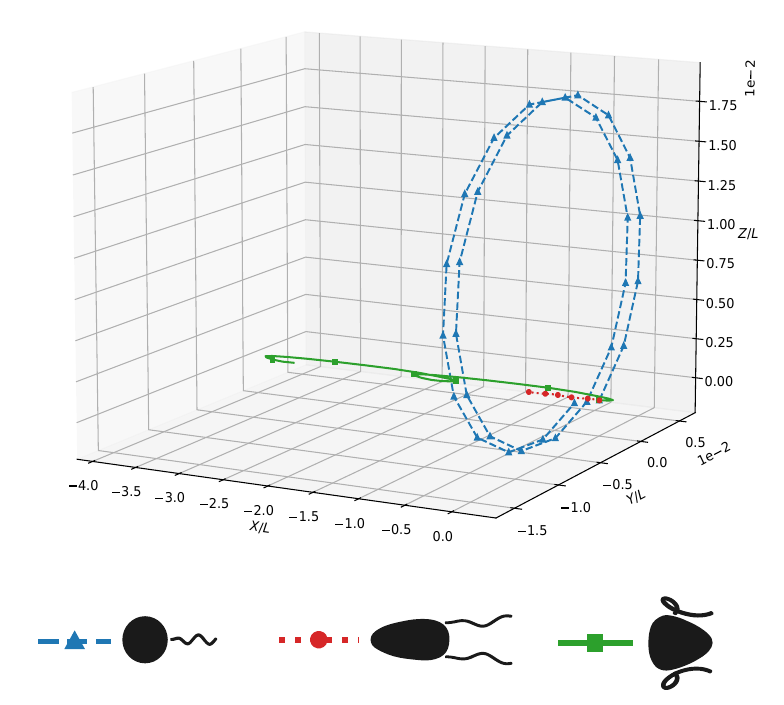}\\
        };
    \end{tikzpicture}
   \caption{Center of mass trajectories of the reference swimmer $S^0$ (blue triangles) compared to optimal swimmers (from \Cref{table:best_models}) optimized for mean velocity (green squares) and mean efficiency (red rounds) for two periods.}
    \label{fig:best_traj}
\end{figure}

\noindent {\bf Swimming trajectory :} As can be seen on \Cref{fig:best_traj}, the optimal swimmers only move in the $e_1$ direction during their stroke, unlike the reference swimmer. This is due to the second flagellum, which is placed symmetrically and gives the swimmer a certain stability. The fastest moving swimmer does not take into account the power dissipated, so we can see the appearance of a particular trajectory where the swimmer moves backwards during the stroke, allowing him to reach a high propulsion speed along $e_1$ while dissipating a lot of power. However, when the power dissipated is taken into account, the swimmer adopts a strategy that allows him to make as few "unnecessary" movements as possible. 

\begin{table}[htpb]
\caption{\label{table:best_models} Comparison between the reference swimmer, $S^0$, introduced in \Cref{table:S0}, and the best swimmers with one and two flagella found through Bayesian optimization using two cost functions \eqref{eq:cost_func}: mean velocity (first column) and mean efficiency (second column). The first row represents the reference swimmer. The second row displays optimized monoflagellated swimmers, and the third row shows optimized biflagellate swimmers.}
\begin{ruledtabular}
  \begin{tabular}{c c c}
     & Mean speed  & Mean efficiency \\ 
  \colrule
       \raisebox{10.5ex}{Reference swimmer $S^0$}  & \multicolumn{2}{c}{\includegraphics[width=5cm]{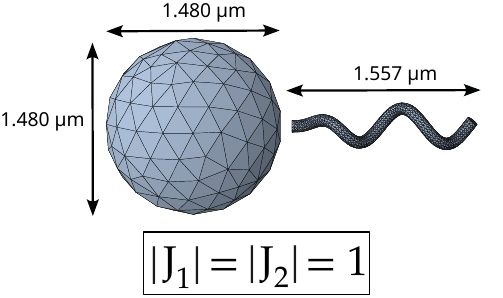}}  \\
  \colrule
     \raisebox{7.5ex}{Best monoflagellated swimmers} &     \includegraphics[width=5cm]{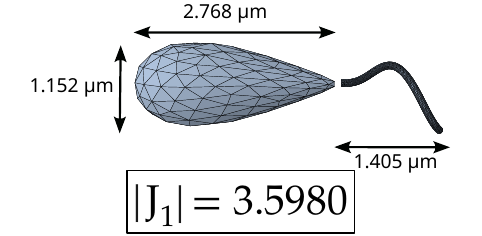} &   \includegraphics[width=5cm]{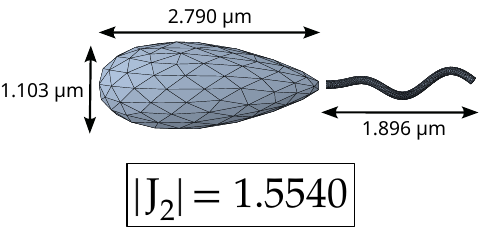} \\
     \raisebox{18ex}{Best biflagellate swimmers}&  \includegraphics[width=3.5cm]{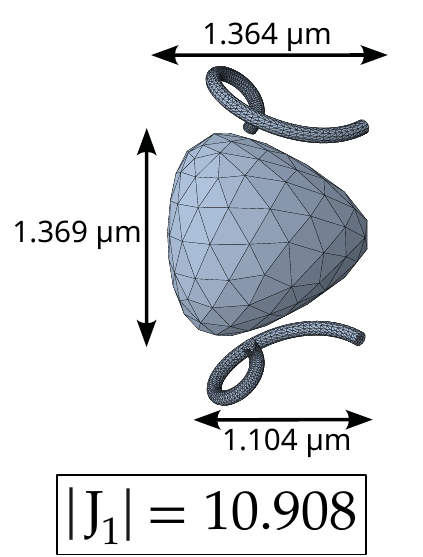} &\raisebox{5ex}{\includegraphics[width=5cm]{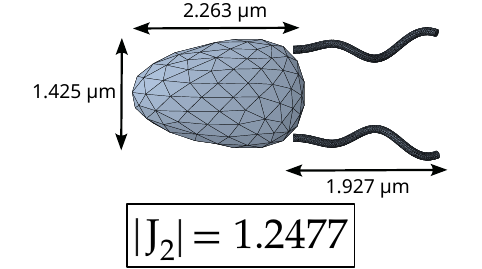}}\\
  \end{tabular}
\end{ruledtabular}
\end{table}

\begin{table}[htpb]
\caption{\label{table:num_res} Swimming and geometrical characteristics of optimized swimmers from \Cref{table:best_models}, normalized relative to the reference swimmer $S^0$ (see \Cref{table:S0_res}).}
\begin{ruledtabular}
  \begin{tabular}{c c c c c c}
   & &  \multicolumn{2}{c}{Optimization of $J_1$} & \multicolumn{2}{c}{Optimization of $J_2$}\\
   \colrule
    & \includegraphics[width=1.5cm]{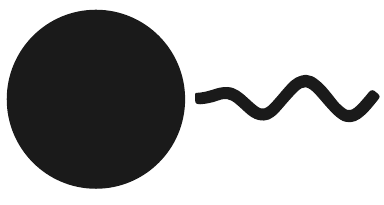} &  \includegraphics[width=1.5cm]{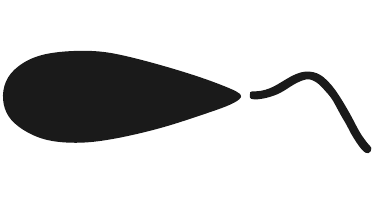} & \includegraphics[width=0.7cm]{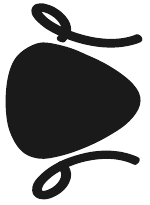} & \includegraphics[width=1.8cm]{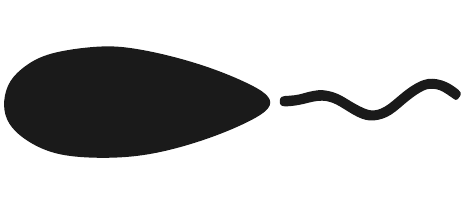} & \includegraphics[width=1.5cm]{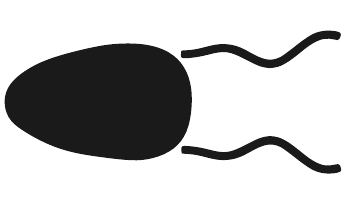}\\
    \colrule
    $\bar{U}_1/\bar{U}^0_1$ & $1$ & $3.5980$ & $10.908$ & $1.5196$ & $2.6046$\\
    $\bar{\Omega}_1/\bar{\Omega}^0_1$ & $1$ & $5.3728$ & $1.7061$ & $1.3660$ & $1.3993$\\
    $\bar{P}/\bar{P}^0$ & $1$ & $3.9366$ & $15.144$ & $0.9785$ & $2.0874$\\
    \colrule
    $\lambda/\lambda^0$ & $1$ & $2.5705$ & $1.0169$ & $1.7162$ & $1.9069$\\
    $R^t/R^{t0}$ & $1$ & $4.9985$ & $4.9965$ & $1.0985$ & $1.1745$\\
    $\alpha/\pi$ & $0$ & $2.5485e{-3}$ & $0.3686$ & $2.5420e{-3}$ & $0.1868$\\
    $\beta/\pi$ & $0$ & $1.7318e{-3}$ & $0.4994$ & $3.0524e{-2}$ & $3.9592e{-2}$\\
    $\gamma/\pi$ & $0$ & $-2.1914e{-3}$ & $0.2518$ & $2.3316e{-3}$ & $1.2893e{-2}$ \\
    $\delta/\pi$ & $0$ & $5.9362e{-3}$ & $-0.2911$ & $4.7056e{-3}$ & $1.5362e{-2}$ \\ 
  \end{tabular}
\end{ruledtabular}
\end{table}

\pagebreak
\section{Discussions}\label{sec:4}

\begin{figure}[htpb]
    \centering
          \includegraphics[width=17cm]{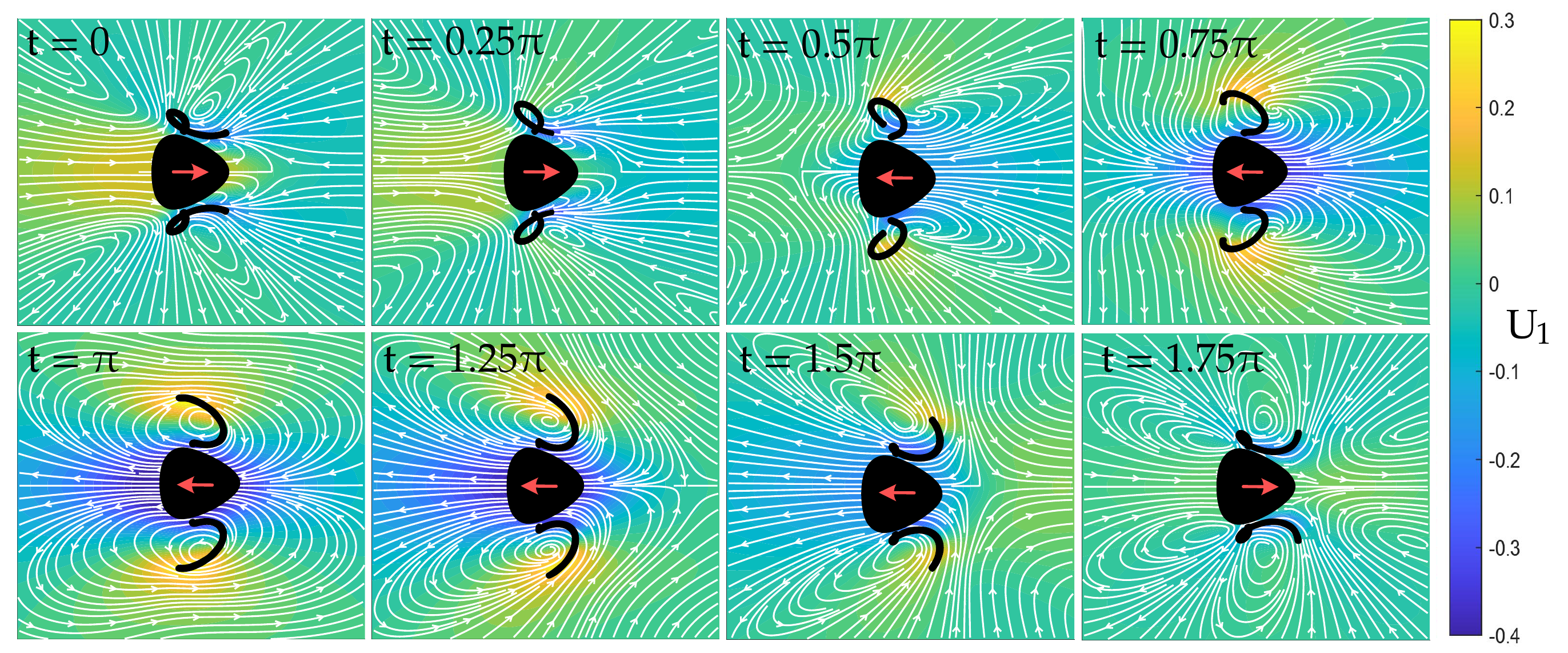}
  \caption{Time sequence of the velocity field evolution throughout the stroke cycle of the bullhead swimmer. The fluid streamlines in the plane $y=0$ are depicted in white. The colormap represents the fluid velocity along $e_1$. The red arrow indicates the direction of the swimmer's displacement. ($t=0$-$0.25\pi$) The swimmer exhibits a small recoil. ($t=0.5\pi$-$1.5\pi$) The swimmer moves forward with a peak velocity at $t=\pi$. ($t=1.75\pi$) The swimmer experiences a small recoil similar to the start of the movement.}
    \label{fig:flow2flagU}
\end{figure}

\begin{figure}[htpb]
    \centering
          \includegraphics[width=17cm]{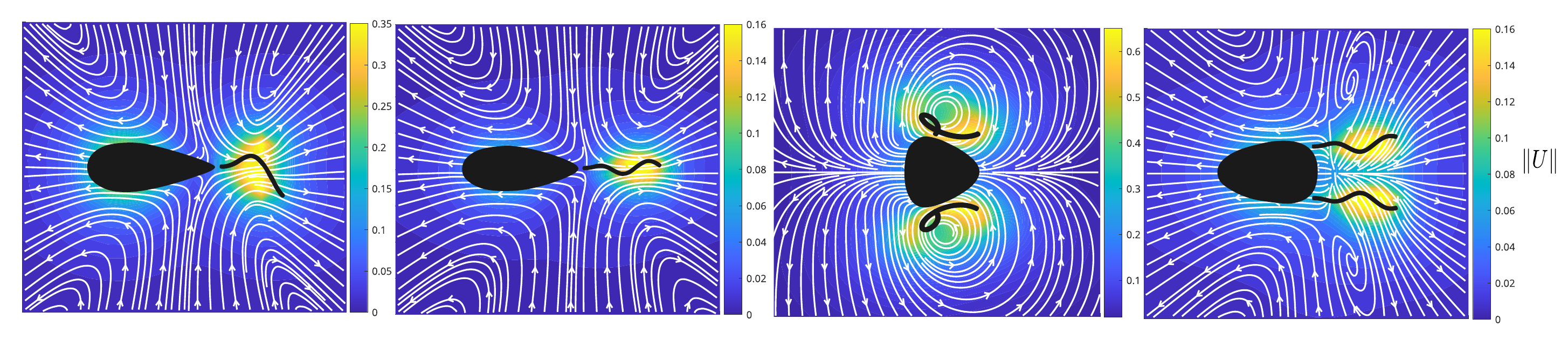}
 \caption{Time-averaged velocity field around optimized swimmers: (first) optimized monoflagellated swimmer for $J_1$, (second) optimized monoflagellated for $J_2$, (third) optimized biflagellate for $J_1$, and (fourth) optimized biflagellate for $J_2$. The fluid streamlines in the plane $y=0$ are depicted in white. The colormap illustrates the magnitude of the fluid velocity, with each representation adjusted to its specific scale due to significant variations in velocity intensities.
 }
    \label{fig:flows}
\end{figure}

Our results on the shape optimization of flagellated microswimmers through Bayesian optimization offer valuable insights into the efficiency morphology at the microscale. A notable finding from our study is the significant influence of flagellar number and positioning on microswimmer head shapes as can be seen between \Cref{fig:best_J1_1flag} and \Cref{fig:best_J1_2flag}, and \Cref{fig:best_J2_1flag} and \Cref{fig:best_J2_2flag}.  However, our study extends beyond existing research by adopting a bigger admissible shape space. By doing so, we discover novel  and better optimal shapes that diverge from conventional ellipsoidal designs \cite{kenta_optisperm, phan-thien_tran-cong_ramia_1987,Lauga_2009, acemoglu2014effects}, highlighting the importance of considering a new range of morphological possibilities.\\

\noindent
An important characteristic of flagella, revealed through optimization, is that while larger flagellar amplitudes result in higher speeds, they also lead to increased power consumption, as detailed in \Cref{table:num_res}. Amplitude and wavelength also play an important role in the swimmer's motion. This counterintuitive phenomenon was observed both experimentally and analytically, where the fastest magnetic helical swimmers in \cite{walker2015optimal} had a flagellum close to one full turn. While there are differences in motility - as this paper considers a torque-free swimmer, in contrast to the torque-driven helical propeller in \cite{walker2015optimal} - the same geometric property can be observed in both cases. This trade-off has significant implications for microswimmer design, emphasizing the need for nuanced optimization strategies that balance velocity with efficiency.\\

\noindent In \Cref{table:num_res}, we note certain phenomena common to the mono- and biflagellated cases. In the $J_1$ case \eqref{eq:cost_func}, the angular velocities and the power dissipated are higher than in the other cases, and the flagellum amplitude is almost equal to the upper limit chosen. In the $J_2$ case, on the other hand, we obtain lower angular velocities and flagellum amplitudes. We can also see that the orientation of the flagella, i.e., $\gamma$ and $\delta$, are very small with the exception for the bullhead swimmer. Additionally, we observe that the optimized monoflagellated swimmer is more efficient than the optimized biflagellate swimmer ($1.5540 > 1.2477$). Indeed, the biflagellate swimmer dissipates slightly more than twice the power of the monoflagellated swimmer. Moreover, the addition of a second flagellum does not result in a doubling of the translational velocity. An important difference can be seen in \Cref{fig:best_traj} between mono- and biflagellate swimmers. Swimmers with a single flagellum produce helical trajectories. This phenomenon, known as the \textit{helix theorem}, is described in more detail in \cite{helixth_cf, helixth_pnas}. The addition of a second symmetrical flagellum results in a helical trajectory with zero amplitude, unlike the monoflagellated swimmer. This characteristic is important to consider in the context of maneuverability. As shown in \cite{noblinelife}, the second flagellum plays a crucial role in achieving great maneuverability, particularly in the reorientation of microorganisms.\\

\noindent In contrast to prior works that often focus on optimizing individual parameters incrementally \cite{acemoglu2014effects, Bet_efficient_shapes_2017}, or even testing different types of flagella \cite{AD2022} our approach considers the simultaneous optimization of complex head and flagella parameters. This methodology allows us to achieve generalized optimization results, as demonstrated by the distinct morphologies observed in single and biflagellated microswimmers. Our findings suggest that by incorporating considerations of power dissipation, we can achieve more streamlined head shapes. These shapes are both flattened and puffed in certain areas, similar to the head of a human spermatozoon \cite{SMITH_2009}, with the difference that the attachment point of our flagella is in the puffed area.\\

\noindent Two aspects are particularly distinctive for the bullhead swimmer: its shape and its trajectory. Unlike other swimmers, a low-amplitude recoil appears during its stroke. A study of the fluid in the neighbourhood of this swimmer at eight different instants during a stroke is described in \Cref{fig:flow2flagU}. We observe a low-amplitude backward movement at $t \in \{0, \frac{\pi}{8}, \frac{7\pi}{8}\}$, while during the rest of the stroke, the forward displacement is much greater. Additionally, \Cref{fig:flow2flagU} closely resembles the time sequence of the velocity field observed for \textit{Chlamydomonas reinhardtii} in Figure 3 of \cite{pusherGuasto}, which exhibits a \textit{puller} behavior. In \Cref{fig:flows}, which presents time-averaged velocity fields, the other swimmers display a typical \textit{pusher} behavior, while the bullhead swimmer demonstrates a behavior closer to \textit{neutral}, as illustrated in Figure 1(a) of \cite{chaithanya_2021}.

\section{Conclusion and perspectives}\label{sec:5}

The implementation of a more general parametric optimization of shape has enabled the discovery of new optimal shapes for microswimmers, which were previously represented almost exclusively by ellipsoidal forms. The obtained characteristics provide a better understanding of the impact of the shapes of the heads and flagella, as well as their number, on microswimming performance.\\

\noindent Future extensions of this work could lead to a better understanding of the forms of microswimmers that can be found in nature or improve the design of microrobots. While our study has provided new results on the optimal shapes of flagellated microswimmers, many avenues for further exploration remain open. For instance, incorporating advanced techniques such as machine learning, or different shape spaces as B-spline in \cite{Guo_Zhu_Liu_Bonnet_Veerapaneni_2021}, could reveal new optimal designs. Additionally, investigating more complex environments, such as swimmers immersed in viscoelastic fluids or accounting for the elasticity of the swimmer, holds significant potential, particularly in microrobotics. This is especially relevant in the medical field, where microrobots must navigate multiple fluids with varying complex characteristics (often non-Newtonian), which can cause significant reductions in swimming accuracy \cite{roboticsLi}. The methods presented in this paper could be instrumental in designing microrobots that maintain their effectiveness in such challenging environments.\\

\noindent {\bf Supplemental Material:} See Supplemental Material (video files) available with the submission on arXiv for videos of the movements associated with the optimized swimmers and a race between the best swimmers and the reference swimmer.

\begin{acknowledgments}
The authors acknowledge the support of the French Agence Nationale de la Recherche
(ANR), under Grant No.
ANR-21-CE45-0013, Project NEMO. 
\end{acknowledgments}

\appendix

\section{Benchmark swimmers}\label{app:littswimmers}
\begin{figure}[H]
    \centering
    \begin{tikzpicture}
        \matrix[matrix of nodes, column sep=2cm]{
          \includegraphics[width=5cm]{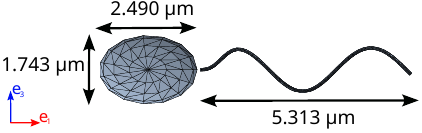} &
          \includegraphics[width=5cm]{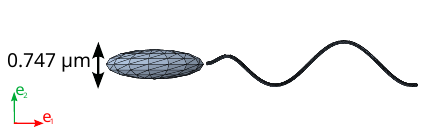}\\
        };
    \end{tikzpicture}
    \caption{Optimal swimmer, in frame $(e_1, e_3)$ (left) and in frame $(e_1, e_2)$ (right), according to the inverse efficiency \eqref{eq:inv_eff} from \cite{phan-thien_tran-cong_ramia_1987} by optimizing only the radius of the ellipsoidal head. The geometric parameters are presented in \Cref{table:87_param}.}
    \label{fig:swimmer_best_87}
\end{figure}

Detailed presentation of swimmers from the literature are depicted in this appendix. The best swimmer, according to the inverse efficiency \eqref{eq:inv_eff}, obtained in \cite{phan-thien_tran-cong_ramia_1987} by optimizing only the parameters related to the head is described in \Cref{table:87_param} and its mesh is displayed in \Cref{fig:swimmer_best_87}. 

\begin{table}[H]
\caption{\label{table:87_param} Geometric parameters describing the head and the flagella of the optimal bacteria in \Cref{fig:swimmer_best_87} from \cite{two_flagella_shum} by optimizing only the radius of the ellipsoidal head.}
\begin{ruledtabular}
  \begin{tabular}{c c c}
     Parameter & Value (Dimensionless) & Value (Dimensional) \\ 
     \colrule
     $R_1$ & $(0.7\times 0.3)^{-1/3}$ & $1.245$ µm\\ 
      $R_2$ & $0.3 \times R_1$ & $0.374$ µm\\
     $R_3$ & $0.7 \times R_1$ & $0.872$ µm\\
    $\bar{A}$& $1$ & $0.74$ µm\\
     $L$ & $10\times  \bar{A}$ & $7.4$ µm \\
     $r$ & $0.05 \times \bar{A}$ & $0.037$ µm\\
      $\lambda$ & $4.7863$ & $3.542$ µm\\
     $R^t$ & $\lambda/(2\pi)$ & $0.564$ µm \\
     $k_E$ & $2\pi/\lambda$ & $1.773$ $\text{µm}^{-1}$\\
  \end{tabular}
\end{ruledtabular}
\end{table}

\begin{figure}[H]
    \centering
    \begin{tikzpicture}
        \matrix[matrix of nodes]{
          \includegraphics[width=6cm]{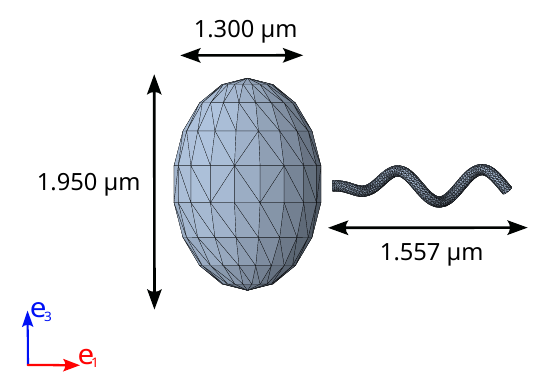} &
          \includegraphics[width=6.4cm]{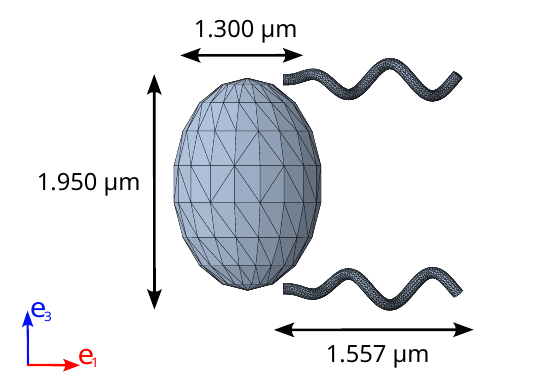}\\
        };
    \end{tikzpicture}
    \caption{(left) Mesh of the reference monoflagellated swimmer from \cite{two_flagella_shum} with $\alpha=0$ and the geometric parameters in \Cref{table:geo_shum}. (right) Example mesh from \cite{two_flagella_shum} for biflagellated swimmer with $\alpha=0.4\pi$ and the geometric parameters in \Cref{table:geo_shum}.}
    \label{fig:swimmer_shum}
\end{figure}

In \cite{two_flagella_shum}, microswimmers are described by a fixed ellipsoidal head, with only the angle $\alpha$ varying to yield the best swimmer according to the average efficiency \eqref{eq:cost_func} at $\alpha = 0.4\pi$. The characteristics of the swimmer's head and flagella are presented in \Cref{table:geo_shum}. \Cref{fig:swimmer_shum} shows the reference swimmer from \cite{two_flagella_shum} on the left, used to normalize speeds and powers, and the best swimmer with an angle $\alpha$ equal to $0.4\pi$ on the right.

\begin{table}[H]
\caption{\label{table:geo_shum} Geometric parameters describing the head and the flagella of the bacteria in \cite{two_flagella_shum}.}
\begin{ruledtabular}
  \begin{tabular}{c c c}
     Parameter & Value (Dimensionless) & Value (Dimensional) \\ 
     \colrule
     $R_1$ & $0.874$ & $0.65$ µm \\ 
     $R_2$ & $0.874$ & $0.65$ µm \\
     $R_3$ & $1.5R_1$ & $0.975$ µm \\
     $L$ & $3.0$ & $2.2$ µm \\
     $r$ & $0.067$ & $0.05$ µm \\
     $R^t$ & $0.2$ & $0.15$ µm \\
     $\lambda$ & $1.0$ & $0.74$ µm \\
     $k_E$ & $0.333\times 2\pi/\lambda=2.09$ & $2.8$ $\text{µm}^{-1}$
 \\
     $l$ & $2r=0.134$ & $0.1$ µm \\
  \end{tabular}
\end{ruledtabular}
\end{table}

\section{Validation of the BEM.}\label{app:valBEM}

\subsection{Monoflagellated microswimmer} 
\begin{figure}[htpb]
    \centering
    \includegraphics[width=12cm]{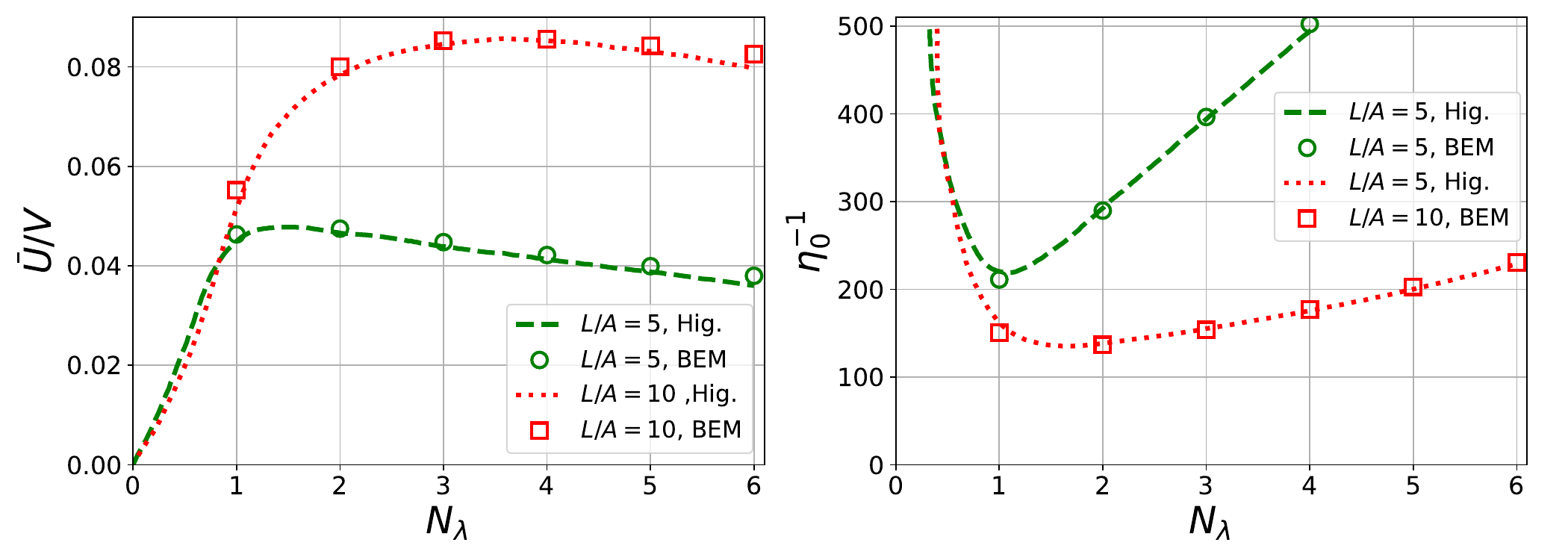}
    \caption{Comparison of normalized mean translational velocity (left) and the inverse efficiency (right) between the results obtained in \cite{higdon_helical_waves} (dashes) and our simulations for $L/A=5$ (round) and $L/A=10$ (square).}
    \label{fig:validation_bem_87}
\end{figure}

We reproduce here a case from the literature \cite{higdon_helical_waves} where Slender-Body-Theory (SBT) was used to solve a self-propulsion problem. The aim is to recover the mean translational velocity and the inverse efficiency for a monoflagellated bacterium according to the number of wavelength $N_{\lambda}$. This number is related to the total tail length $L$ by
\begin{equation*}
    L=\int_{0}^{\lambda N_{\lambda}}\left\|\left(\frac{dx(s)}{ds}, \frac{dy(s)}{ds}, \frac{dz(s)}{ds}\right)\right\|_{2} ds ,
\end{equation*}
where $(x(s), y(s), z(s))$ are given by \eqref{eq:flag}. The mean translational velocity is normalized by $V=\omega/k_E$ \eqref{eq:flag}-\eqref{eq:stokes} and the inverse efficiency is given, here, by \eqref{eq:inv_eff}. The head of the swimmer is taken as sphere of unit radius $A=1$. The geometry of the swimmer considered is referenced in \Cref{table:geo_pt} where $\alpha=\gamma=\beta=\delta=0$. \Cref{fig:validation_bem_87} shows the comparison between the different results obtained for $L/A=5$ and $L/A=10$ ratio, and those in the paper \cite{higdon_helical_waves}.

\begin{table}[htpb]
\caption{ \label{table:geo_pt} Geometric parameters describing the flagellum of the bacteria in \cite{higdon_helical_waves}.}
\begin{ruledtabular}
  \begin{tabular}{c c c}
     Parameter & Value (Dimensionless) & Value (Dimensional)\\ 
     \colrule
     $A$ & $1$ & $0.74$ µm\\
     $r$ & $0.02A$ & $0.0148$ µm \\
     $l$ & $2r$ & $0.0296$  µm\\
     $R^t$ & $\lambda/(2\pi)$ &\\
     $k_E$ & $2\pi/\lambda$ & \\
  \end{tabular}
\end{ruledtabular}
\end{table}

\subsection{Biflagellated microswimmer.}
\begin{figure}[htpb]
    \centering    \includegraphics[width=11cm]{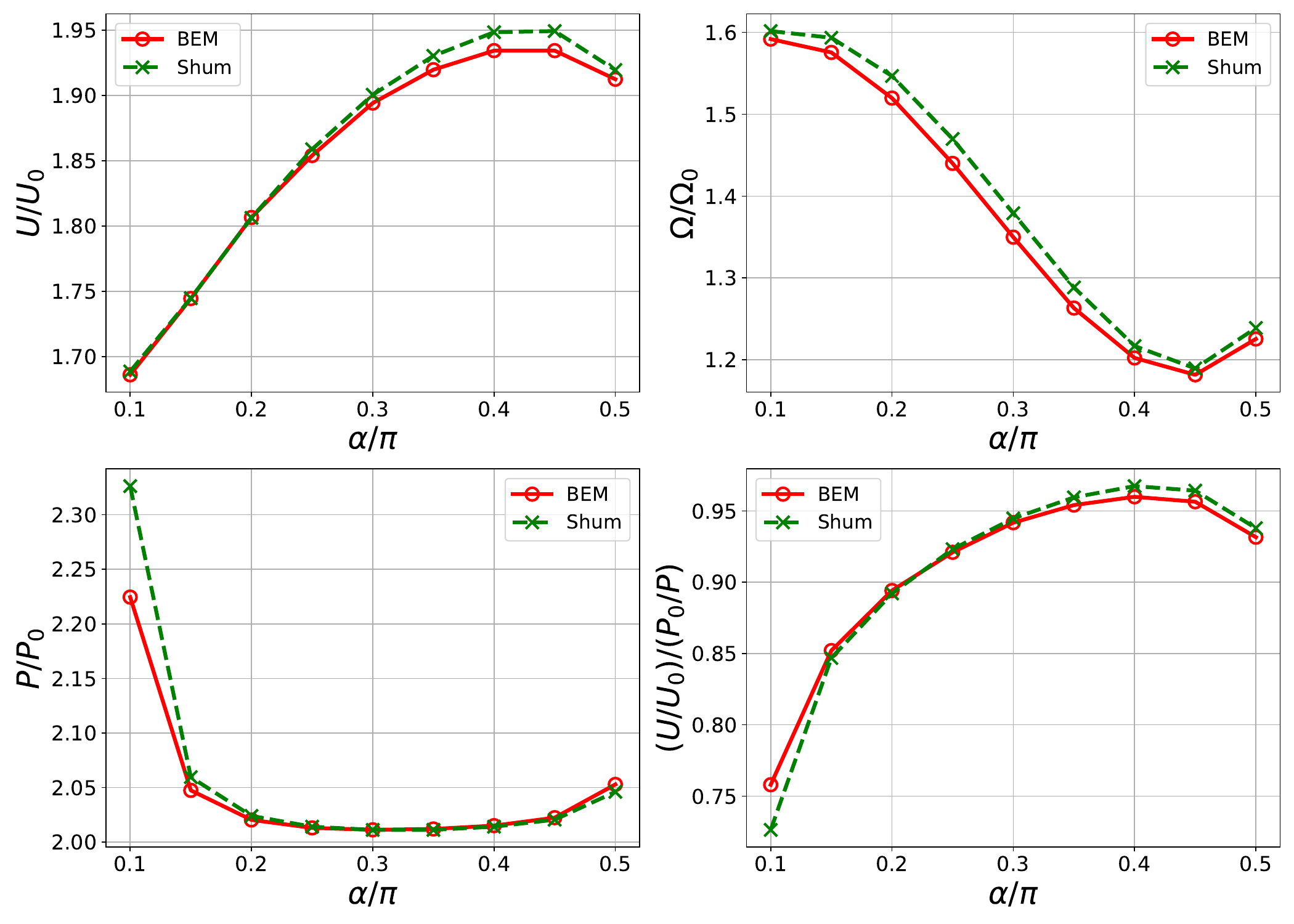}
    \caption{Comparison of normalized mean translational velocity (top left), normalized mean rotational velocity (top right), normalized mean power dissipation (bottom left) and normalized mean efficiency (bottom right) between the results obtained in \cite{two_flagella_shum} (cross), and our simulations (round).}
    \label{fig:validation_bem_shum}
\end{figure}

We reproduce a case from the literature \cite{two_flagella_shum} where the BEM was used to solve a self-propulsion problem. The aim is to recover the mean translational and rotational velocities, mean power dissipation, and efficiency for a biflagellate bacterium, normalizing these quantities by those obtained for a monoflagellated bacterium with the same head and flagellum shape, as presented on the left in \Cref{fig:swimmer_shum}.
 The different quantities are described by equations \eqref{eq:mean_P} and \eqref{eq:mean_vel}, with respect to the new reference swimmer chosen in \Cref{table:geo_shum}.
\Cref{fig:validation_bem_shum} shows the comparison between the results and those in the paper \cite{two_flagella_shum}. The shape of the head is defined by \eqref{eq:ellipsoidal_head}. In addition, the second flagellum is considered to be created by a rotation of $\pi$ around the axis of propulsion $e_1$ of the first flagellum, whose junction point is determined solely by the angle $\alpha$. The geometric parameters used for the simulations are shown in \Cref{table:geo_shum} where $\gamma=\beta=\delta=0$.

\newpage

\section{Shape Variability with FFD}\label{app:ffdillu}

\begin{figure}[htpb]
    \centering
    \includegraphics[width=15cm]{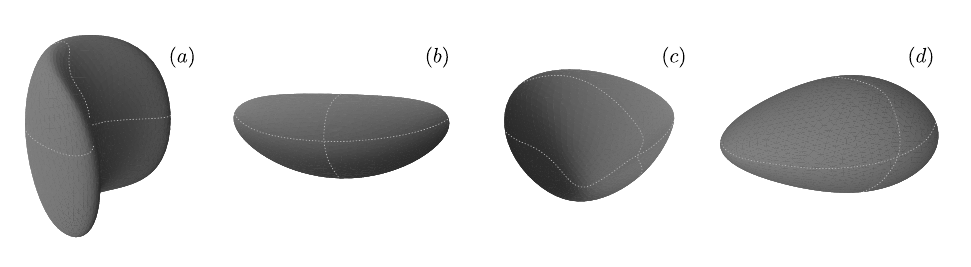}
    \caption{Examples of random shapes generated using FFD by applying displacements to control points within the set $\mathcal{V}$. Shapes $(a)$ and $(c)$ are nonconvex, while shapes $(b)$ and $(d)$ are convex.}
    \label{fig:random_shapes_ffd}
\end{figure}

To illustrate the diversity of shapes belonging to the set $\mathcal{V}$ (defined in equation \eqref{eq:head_ad}) that can be generated using the FFD method, we present a collection of random shapes with varying characteristics in \Cref{fig:random_shapes_ffd}. These shapes are obtained by displacing control points, with same characteristics as in \Cref{app:Numimpl}, while respecting the constraints of the set $ \mathcal{V}$. The geometric constraints, especially symmetry, play a significant role in shaping the outcomes. Despite these constraints, \Cref{fig:random_shapes_ffd} demonstrates the ability to produce a wide range of shapes, including both nonconvex (e.g., shapes $(a)$ and $(c)$) and convex shapes (e.g., shapes $(b)$ and $(d)$). Furthermore, the four examples exhibit distinct curvatures and geometrical features, highlighting the flexibility of the FFD method in generating diverse and complex forms. Additional shapes are presented in \Cref{app:bestswimmer}, where the different optimized head shapes resulting from the FFD method are displayed.

\section{Numerical Implementation}\label{app:Numimpl}

\noindent{\bf Stroke discretization :} In the simulations, the number of BEM evaluations to compute the different mean values \eqref{eq:mean_P} and \eqref{eq:mean_vel} is set to four, obtained by rotating the flagella around their axis by angles of $\frac{\pi}{2}t$ with $t\in\{0,1,2,3\}$. Choosing $6$, $8$, or $12$ evaluations does not significantly change the average values (relative error less than $10^{-3}$), thus $4$ evaluations seems be a good choice especially considering the high computational cost of BEM coupled with the optimization routine.\\

\noindent\textbf{Volume Computation:} To compute the volume of different shapes, we use a boundary integral method. According to the divergence theorem:
\begin{equation*}
\int_{S} (\nabla \cdot F) \, dV = \int_{\partial S} F \cdot n \, dS,
\end{equation*}
where $ F $ is a vector field and $ n $ is the normal vector on the surface $ \partial S $. By choosing $ F(x,y,z) = \frac{1}{3} (x,y,z) $, we obtain:
\begin{equation*}
|S| = \int_{S} dV = \int_{\partial S} \frac{1}{3} (x, y,z) \cdot n \, dS.
\end{equation*}
Since we work with triangular element meshes, we approximate the volume as:
\begin{equation*}
|S| \approx \frac{1}{3} \sum_{\text{triangle}} \int_{\text{triangle}} (x,y,z) \cdot n \, dS.
\end{equation*}

\noindent{\bf Symmetry constraint :} To obtain a varied shape space for the heads of microswimmers while limiting the dimension of the problem due to computational costs and optimization performance, we consider four control points along the three axes, i.e., $M=M_1\times M_2\times M_3=64$ with $M_1=M_2=M_3=4$. As we are interested in the movement of the swimmers along the $e_1$ axis, it seems natural to consider symmetries along the $(e_1,e_3)$ and $(e_1,e_2)$ planes. To do this, we apply this symmetry constraint to our control points (since the reference shape we are deforming, $H^0$, itself respects these symmetries). All that remains is to determine the position of $(64-8)/2=14$ control points (as explained above, we're not touching the internal control points). We are working in dimension $d=3$, so the head optimization problem is of dimension $14\times 3=42$. The overall dimension of the optimization problem is $d = 42+p$ where $p$ is the number of flagellum parameters being optimized.\\

\noindent{\bf Trajectory constraint :} Additional constraints must be added to ensure that swimmers, on average, do not move along $e_2$ or $e_3$ as they swim, described by:
\begin{equation}\label{eq:e1constaint}
    |\bar{U}_2|,|\bar{U}_3|\leq \varepsilon_{\text{tol}},
\end{equation}
where $\varepsilon_{\text{tol}}>0$ is a chosen tolerance.\\ 

\noindent{\bf Mesh collisions :} On top of that, a new constraint must be considered to avoid mesh collisions. Depending on the shape of the head and the parameters of the flagella, such collisions may occur. To avoid this, we introduce a new constraint to the optimization problem. The natural constraint would be to consider a binary constraint informing whether or not the parameters are feasible or infeasible. This type of constraint can be complicated to deal with, particularly for complicated problems because of its binary
discrete nature. It is preferable to work with continuous constraints. Here, we use a Monte Carlo (MC) method to quantify this: we draw $N_{\text{MC}}=10^{6}$ points inside a hypercube and look at the number of points that are inside, $N_{\text{in}}$, at least two different convex envelopes associated with the meshes. Considering the convex envelopes of the flagella allows us to check whether there is an intersection that can occur when the tail rotates around its axis. The constraint is defined as:
\begin{equation}\label{eq:mc}
   \frac{N_{\text{in}}}{N_{\text{MC}}}\leq 0.
\end{equation}

\noindent The numerical parameters of the optimization problem are described in \Cref{table:param_opti}. The optimization problem \eqref{eq:cost_func}  becomes,
for $i \in \{1,2\}$,
\begin{equation}\label{eq:cost_func_2}
    \inf_{\substack{S \in \mathcal{S}_{ad}\\N_{in}/N_{MC}\leq 0\\|\bar{U}_2|,|\bar{U}_3|\leq \varepsilon_{tol}}} J_i(S).
\end{equation}

\begin{table}[H]
\caption{ \label{table:param_opti} Constraints optimization parameters.}
\begin{ruledtabular}
  \begin{tabular}{l l l}
     Parameter & Value & Description  \\ 
     \colrule
     $\varepsilon$ & $0.01$ & Tolerance on head volume conservation \eqref{eq:head_ad}. \\
      $\varepsilon_{\text{tol}}$ & $0.001$ & Tolerance on the average trajectory according to $e_1$ \eqref{eq:e1constaint}.\\
     $N_{\text{MC}}$ & $5e{6}$ & Number of points for the MC method to avoid mesh collisions \eqref{eq:mc}.\\
       $M_1,M_2,M_3$ & $4$ & Number of control points by direction $e_1$, $e_2$ and $e_3$.\\
       $\tilde{M}$ & $\Pi_{i=1}^3 M_i-\Pi_{i=1}^3(M_i-2)$ & Number of movable control points, i.e. those at the edge.\\
       $r_i$ & $\frac{4}{2(M_i-1)}-\frac{0.4}{M_i}$ & Bounds in the direction $e_i$ for the control points.\\
       $\mathcal{V}$ & $\Pi_{k=1}^{\tilde{M}}\left(\Pi_{i=1}^3 [-r_i,r_i]\right)$ & Closed set for the control points to prevent mesh collapse in \eqref{eq:head_ad}.\\
       $p$ & & Number of flagellum parameters being optimized.\\
       $d$ & $42+p$ & Dimension of the optimization problem.\\
     $r^{\text{SCBO}}$  & $\min(5000, \max(2000,200\times d))$ & Candidates generated in the trust region of SCBO algorithm.\\ 
     $q^{\text{SCBO}}$ & $15$ & Number of batches of SCBO algorithm.\\
      $L_{\text{init}}^{\text{SCBO}}$ & $1.6$ & Initial length of the trust region of SCBO algorithm.\\
    $L_{\text{min}}^{\text{SCBO}}$ & $0.5^{7}$ & Minimal length of the trust region of SCBO algorithm.\\
   $L_{\text{max}}^{\text{SCBO}}$ & $1.6$ &  Maximal length of the trust region of SCBO algorithm.\\
    $n_{\text{init}}^{\text{SCBO}}$ & $3\times d$ & Number of initial points of SCBO algorithm.\\
    $\tau_s^{\text{SCBO}}$ & $\max(3, \lceil d/10\rceil)$& Sucess rate of SCBO algorithm.\\
    $\tau_f^{\text{SCBO}}$ & $\lceil d/q \rceil$ & Failure rate of SCBO algorithm.\\
  \end{tabular}
\end{ruledtabular}
\end{table}

\section{A collection of the best swimmers}\label{app:bestswimmer}
A set of images illustrating in greater detail the swimmers optimized in this paper in different orientations, as well as their dimensions.

\begin{figure}[H]
    \centering
    \begin{tikzpicture}
        \matrix[matrix of nodes, column sep=2cm]{\includegraphics[width=5cm]{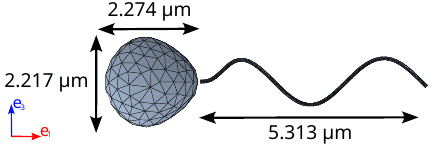} &
        \includegraphics[width=5cm]{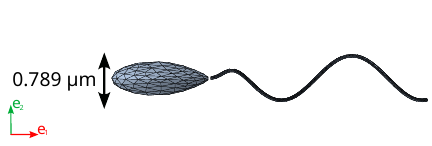}\\          
        };
    \end{tikzpicture}
    \caption{Optimal swimmer obtained by the head optimization of the literature swimmer from \cite{phan-thien_tran-cong_ramia_1987} in the frame $(e_1, e_3)$ (left) and in the frame $(e_1, e_2)$ (right).}
    \label{fig:best_opti_87}
\end{figure}

\begin{figure}[H]
    \centering
    \begin{tikzpicture}
        \matrix[matrix of nodes, column sep=2cm]{
         \includegraphics[width=4cm]{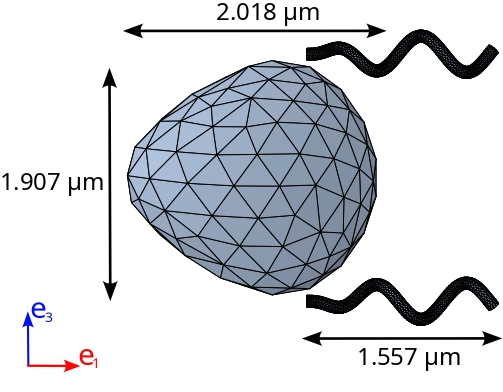} &
        \includegraphics[width=4cm]{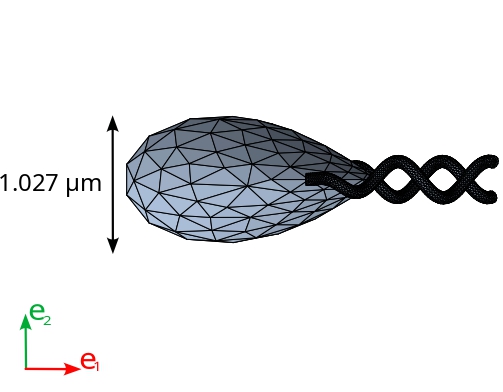}\\
        };
    \end{tikzpicture}
    \caption{Optimal swimmer obtained by the head optimization of the literature swimmer from \cite{two_flagella_shum} in the frame $(e_1, e_3)$ (left) and in the frame $(e_1, e_2)$ (right).}
    \label{fig:best_04pi}
\end{figure}

\begin{figure}[H]
    \centering
    \begin{tikzpicture}
        \matrix[matrix of nodes, column sep=2cm]{
          \includegraphics[width=5cm]{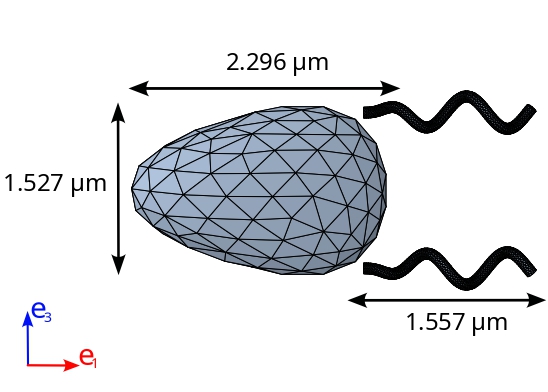} &
          \includegraphics[width=5cm]{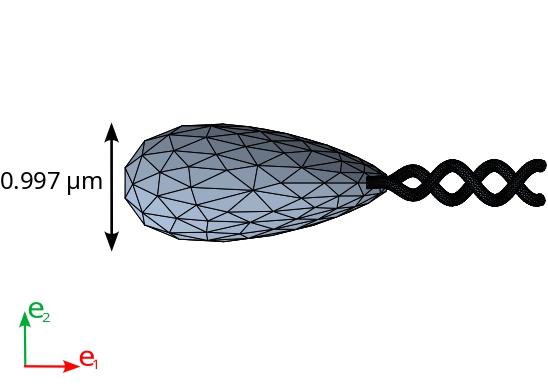}\\  
        };
    \end{tikzpicture}
    \caption{Optimal swimmer obtained by optimising both the head and the $\alpha$ parameter in relation to the literature swimmer from \cite{two_flagella_shum} in the frame $(e_1, e_3)$ (left) and in the frame $(e_1, e_2)$ (right). }
    \label{fig:best_shum}
\end{figure}

\begin{figure}[H]
    \centering
    \begin{tikzpicture}
        \matrix[matrix of nodes, column sep=2cm]{
          \includegraphics[width=5cm]{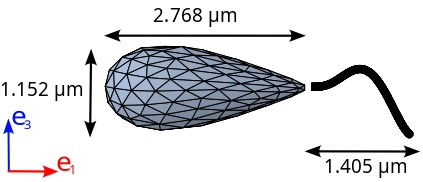} &
          \includegraphics[width=5cm]{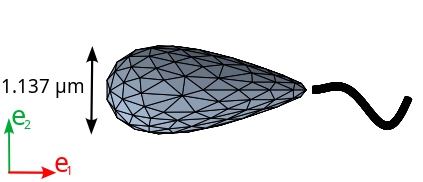} \\
        };
    \end{tikzpicture}
    \caption{Optimal monoflagellated swimmer (\textit{water-drop swimmer}) obtained by optimizing all parameters describing the head and flagellum, taking into account the mean velocity problem $J_1$ \eqref{eq:cost_func} in the frame $(e_1, e_3)$ (left) and in the frame $(e_1, e_2)$ (right).}
    \label{fig:best_J1_1flag}
\end{figure}

\begin{figure}[H]
    \centering
    \begin{tikzpicture}
        \matrix[matrix of nodes, column sep=2cm]{
          \includegraphics[width=5cm]{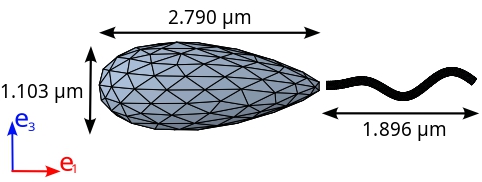} &
          \includegraphics[width=5cm]{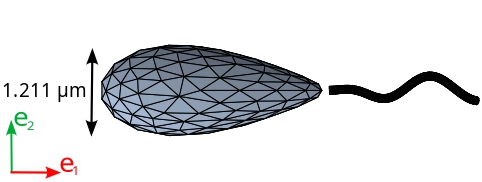}\\
        };
    \end{tikzpicture}
   \caption{Optimal monoflagellated swimmer (\textit{water-drop swimmer}) obtained by optimizing all parameters describing the head and flagellum, taking into account the mean efficiency problem $J_2$ \eqref{eq:cost_func} in the frame $(e_1, e_3)$ (left) and in the frame $(e_1, e_2)$ (right).}
    \label{fig:best_J2_1flag}
\end{figure}

\begin{figure}[H]
    \centering
    \begin{tikzpicture}
        \matrix[matrix of nodes, column sep=2cm]{
          \includegraphics[width=4cm]{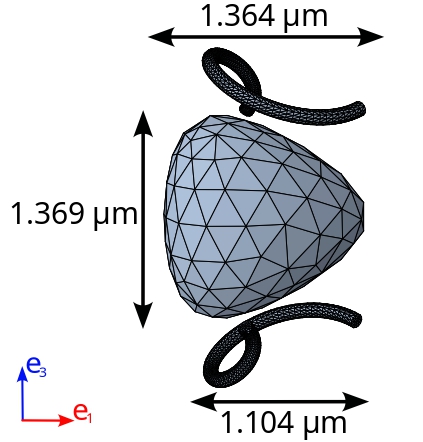} &
         \includegraphics[width=4cm]{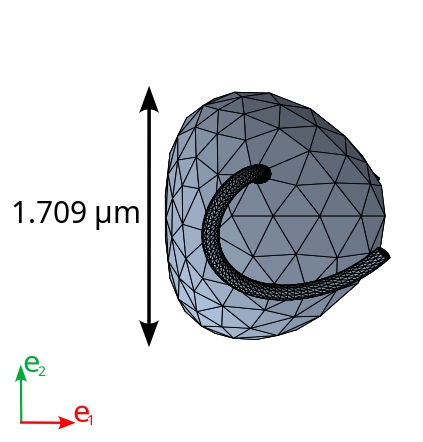}\\
        };
    \end{tikzpicture}
    \caption{Optimal biflagellate swimmer (\textit{bullhead swimmer}) obtained by optimizing all parameters describing the head and flagellum, taking into account the mean velocity problem $J_1$ \eqref{eq:cost_func} in the frame $(e_1, e_3)$ (left) and in the frame $(e_1, e_2)$ (right).}
    \label{fig:best_J1_2flag}
\end{figure}

\begin{figure}[H]
    \centering
    \begin{tikzpicture}
        \matrix[matrix of nodes, column sep=2cm]{
          \includegraphics[width=5cm]{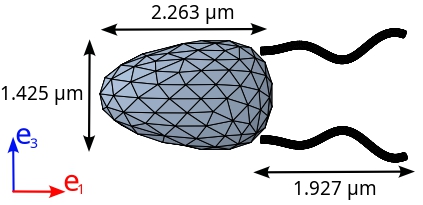} &
          \includegraphics[width=5cm]{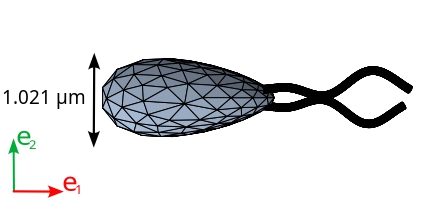}\\
        };
    \end{tikzpicture}
    \caption{Optimal biflagellate swimmer obtained by optimizing all parameters describing the head and flagellum, taking into account the mean efficiency problem $J_2$ \eqref{eq:cost_func} in the frame $(e_1, e_3)$ (left) and in the frame $(e_1, e_2)$ (right).}
    \label{fig:best_J2_2flag}
\end{figure}

\section{Convergence graph}\label{app:cvgraph}
Plots of the different convergence graphs of the SCBO Bayesian optimization algorithm applied to the shape optimization of the various microswimming problems studied in this paper.

\begin{figure}[H]
    \centering
    \begin{tikzpicture}
        \matrix[matrix of nodes]{
          \includegraphics[trim=0.3cm 0.0cm 0.3cm 0.0cm, clip,width=7cm]{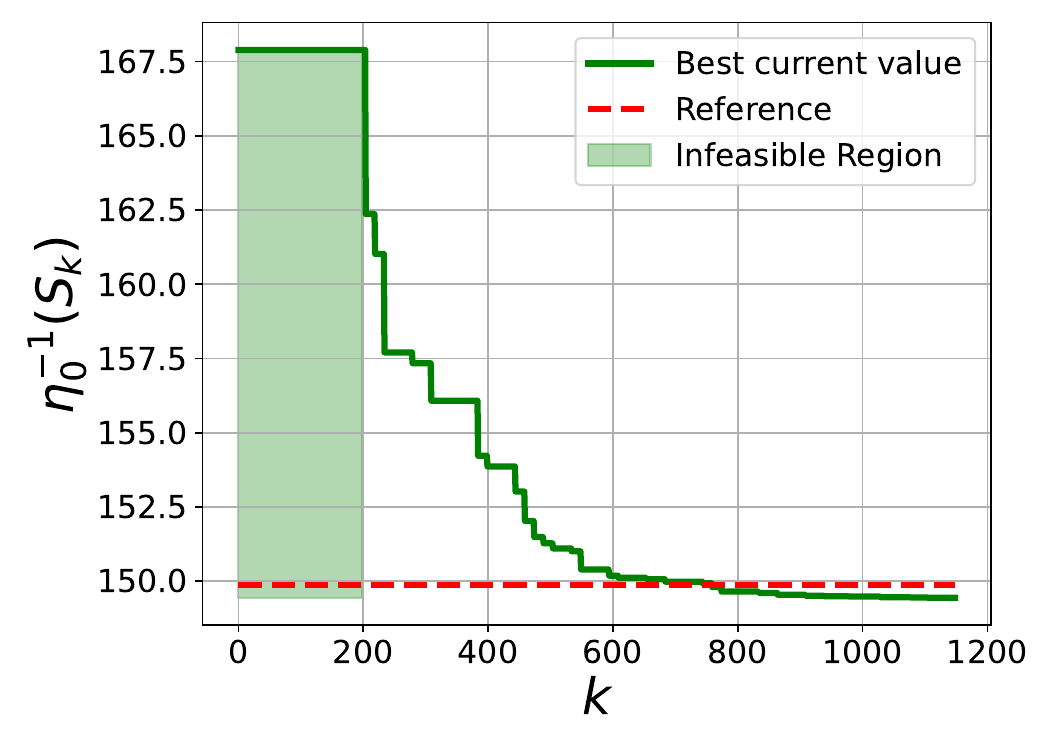} &
          \includegraphics[trim=0.3cm 0.0cm 0.3cm 0.0cm, clip,width=7cm]{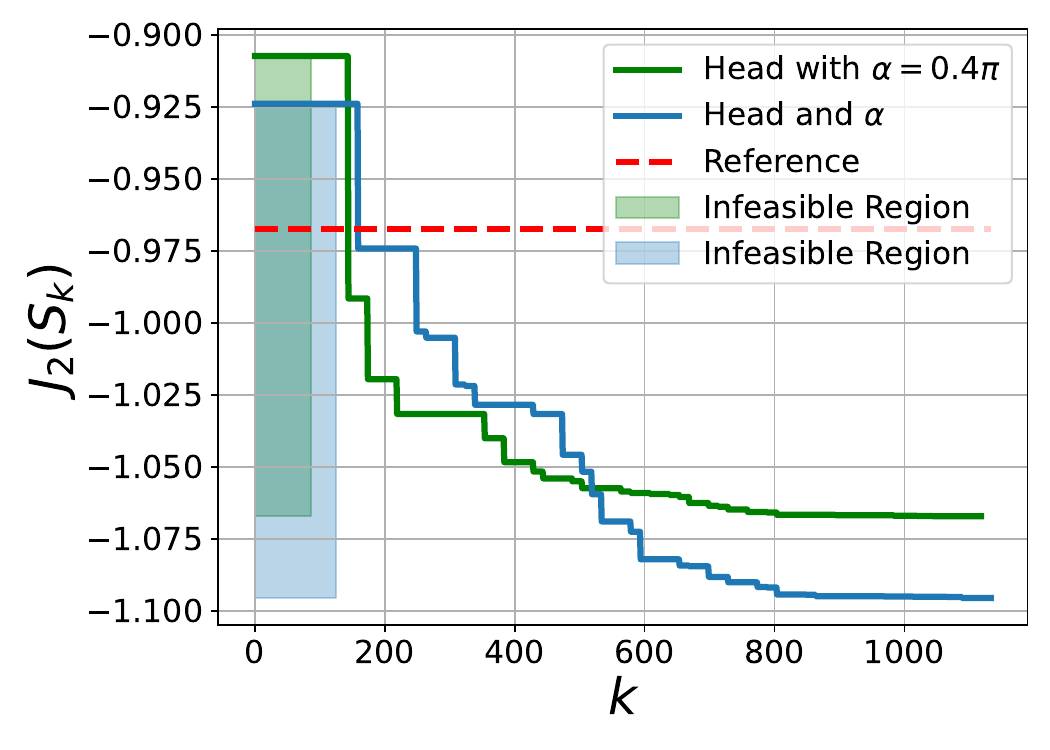}\\
        };
    \end{tikzpicture}
    \caption{The curves are associated with the best current point during the various iterations and batches. The infeasible region corresponds to iterations where constraints have not yet been met. (left) Comparison between the best swimmer in \cite{phan-thien_tran-cong_ramia_1987} (red) and the optimal shape while optimizing the head for the inverse efficiency \eqref{eq:inv_eff}. (right) Comparison between the best swimmer in \cite{two_flagella_shum} (red), the optimal shape while optimizing the head with $\alpha=0.4\pi$ (green), and the optimal shape while optimizing the head and $\alpha$ (blue) for the mean efficiency \eqref{eq:cost_func}. }
    \label{fig:validation_opti}
\end{figure}

\begin{figure}[H]
    \centering
    \begin{tikzpicture}
        \matrix[matrix of nodes]{
          \includegraphics[trim=0.3cm 0.0cm 0.3cm 0.0cm, clip,width=7cm]{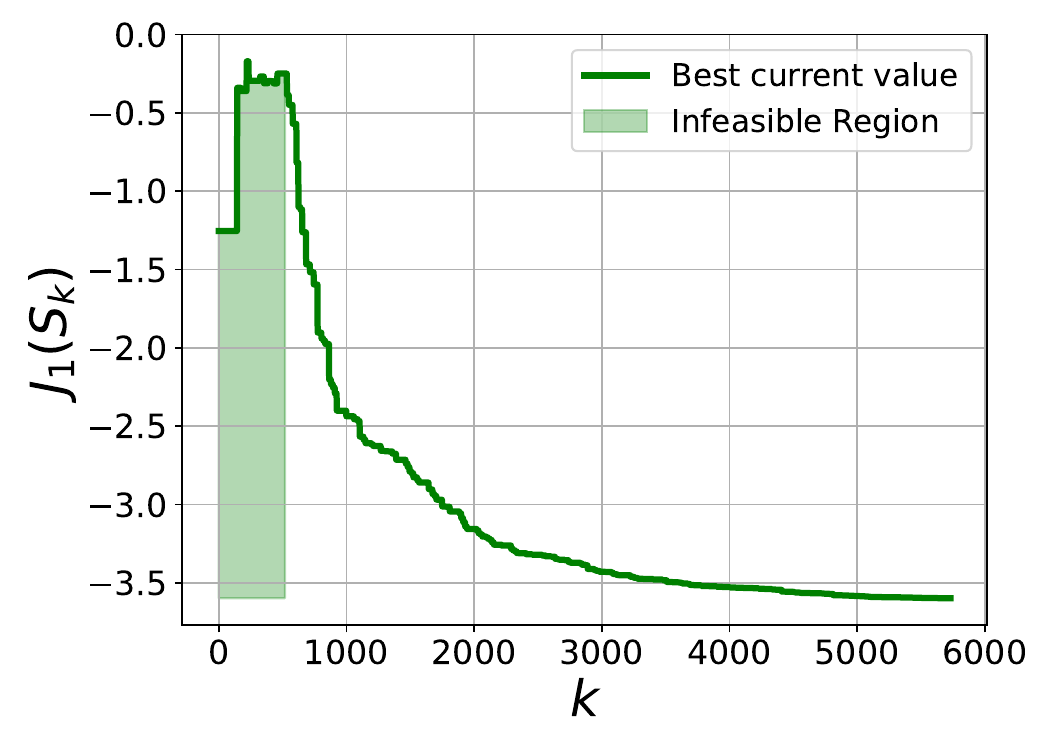} &
          \includegraphics[trim=0.3cm 0.0cm 0.3cm 0.0cm, clip,width=7cm]{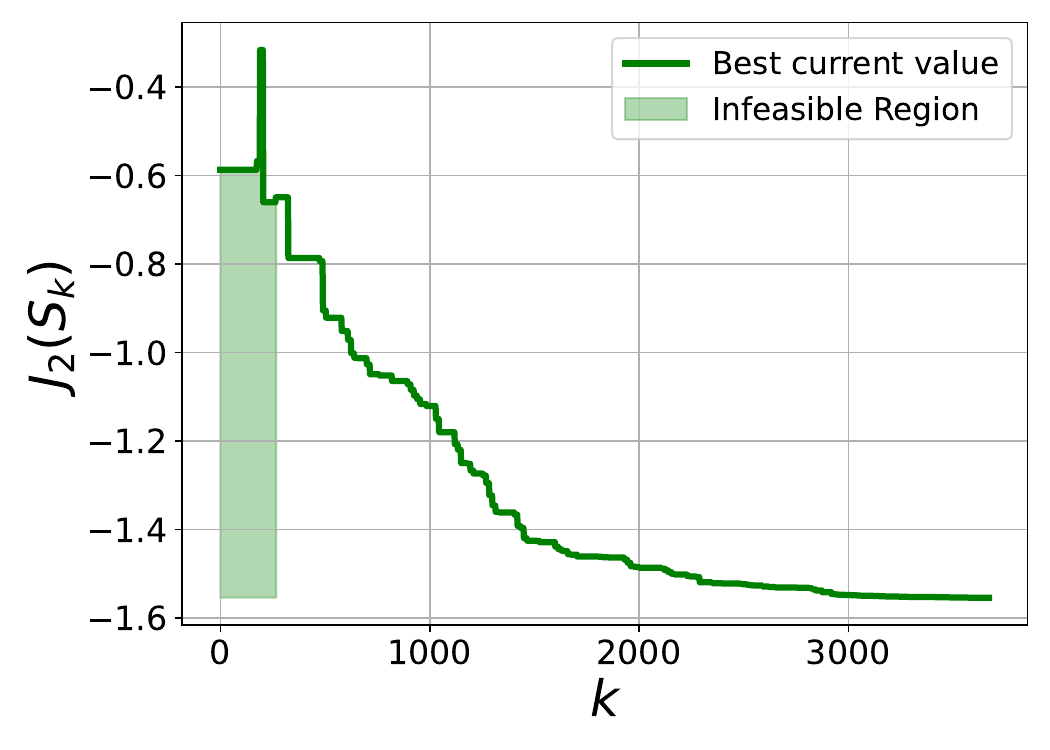}\\
        };
    \end{tikzpicture}
    \caption{The curves are associated with the best current point during the various iterations and batches. The infeasible region corresponds to iterations where constraints have not yet been met. (left) Convergence graph for mean velocity problem $J_1$, \eqref{eq:cost_func}, with one flagellum. (right) Convergence graph for mean efficiency problem $J_2$, \eqref{eq:cost_func}, with one flagellum.}
    \label{fig:best_1flag_cv}
\end{figure}

\begin{figure}[H]
    \centering
    \begin{tikzpicture}
        \matrix[matrix of nodes]{
          \includegraphics[trim=0.3cm 0.0cm 0.3cm 0.0cm, clip,width=7cm]{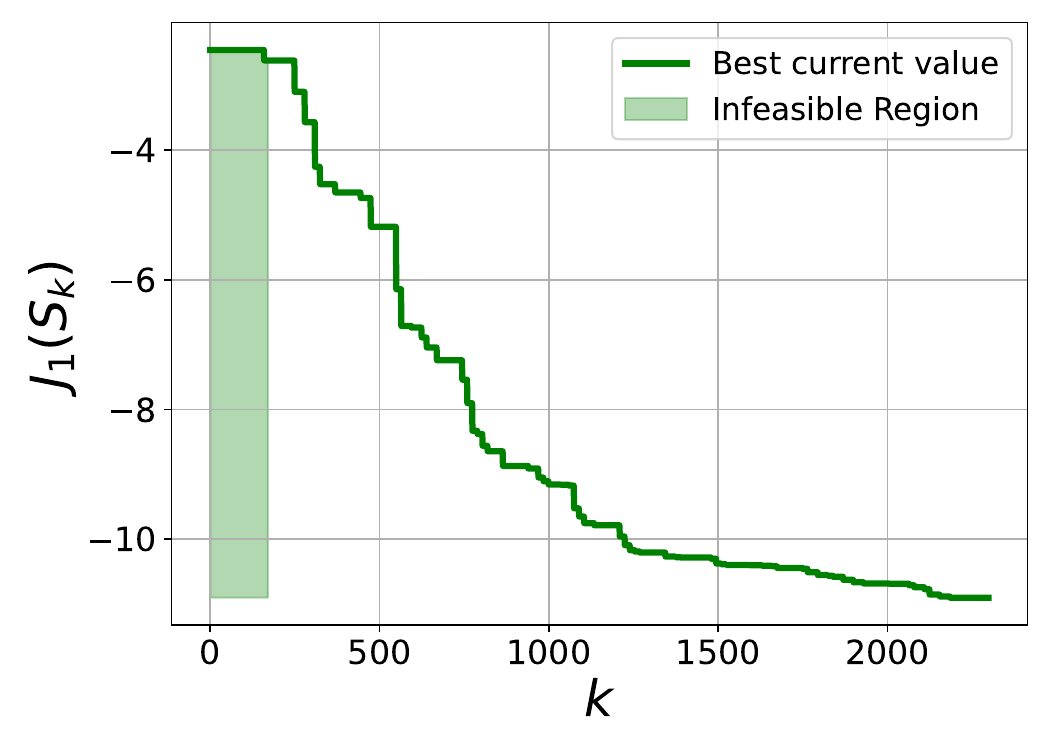} &
          \includegraphics[trim=0.3cm 0.0cm 0.3cm 0.0cm, clip,width=7cm]{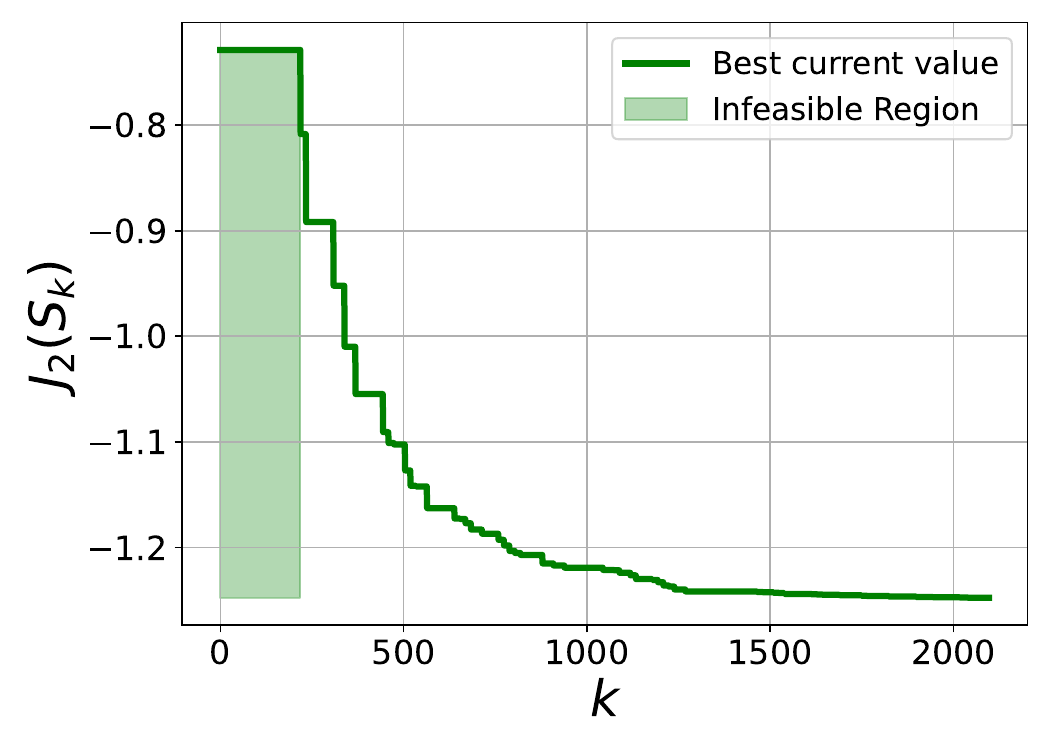}\\
        };
    \end{tikzpicture}
    \caption{The curves are associated with the best current point during the various iterations and batches. The infeasible region corresponds to iterations where constraints have not yet been met. (left) Convergence graph for mean velocity problem $J_1$, \eqref{eq:cost_func}, with two flagella. (right) Convergence graph for mean efficiency problem $J_2$, \eqref{eq:cost_func}, with two flagella.}
    \label{fig:best_2flag_cv}
\end{figure}

\bibliography{biblio.bib}

\end{document}